%% file: ms.tex
\RequirePackage{fix-cm}
\documentclass[twoside,openright,a4paper,11pt,numbers=noenddot,headinclude,footinclude=true,cleardoublepage=empty]{scrreprt}

\input{config.tex}
\usepackage[diplom]{tugrazthesis}
\thesisauthor[Georg Köstenberger]{Georg Köstenberger, BSc}
\thesistitle[Weingarten Calculus]{Weingarten Calculus}
\thesisdate[ ]{November 2020}
\supervisortitle{\germanenglish{Betreuer}{Betreuer}}
\supervisor{%
  Franz Lehner, Assoc.Prof. Dipl.-Ing. Dr.\\

  Institut für Diskrete Mathematik\\

}
\academicdegree{Diplom-Ingenieur}

\curriculum{\germanenglish{ Masterstudium}{Mathematical Data Science}}

\begin{document}
\frenchspacing
\raggedbottom

\maketitle


\cleardoublepage
\pagenumbering{roman}
\setcounter{page}{1}


\tableofcontents

\cleardoublepage

\subfile{frontbackmatter/abstract.tex}

\cleardoubleemptypage
\newpage
\pagenumbering{arabic}
\setcounter{page}{1}

\subfile{chapters/intro/intro.tex}

\subfile{chapters/intro/prelim.tex}

\subfile{chapters/main/main.tex}

\cleardoublepage
\printbibliography
\end{document}

%% file: config.tex
\usepackage[main=english, ngerman]{babel}
\usepackage[T1]{fontenc}

\usepackage{graphicx}
\usepackage{caption}
\graphicspath{{images/}{../images/}{../../images/}} 

\usepackage{booktabs}

\usepackage{amsmath,amssymb,amstext,amsthm,bm,dsfont}

\newtheoremstyle{theorem}
  {\topsep}
  {\topsep}
  {\itshape}
  {0em}
  {\scshape}
  {. \\}
  { }
  {\thmname{#1}\thmnumber{ #2}\thmnote{ (#3)}}

  \theoremstyle{theorem}
  \newtheorem{theorem}{Theorem}
  \newtheorem{lemma}[theorem]{Lemma}
  \newtheorem{corollary}[theorem]{Corollary}
  \newtheorem{mprop}[theorem]{Proposition}
  
  \theoremstyle{definition}
  \newtheorem{mdef}[theorem]{Definition}
  \newtheorem*{remark}{Remark}
  \newtheorem{example}[theorem]{Example}

  \newcommand{\mb}[1]{ \ensuremath{\mathbb{#1}} }
  
  \newcommand{\set}[1]{\left\{#1\right\}}
  \newcommand{\multiset}[1]{\left\{\!\!\left\{#1\right\}\!\!\right\}}
  \newcommand{\mker}[1]{\mathrm{ker}\left(#1\right)}
  \newcommand{\mim}[1]{\mathrm{im}\left(#1\right)}
  \newcommand{\mdim}[1]{\mathrm{dim}\left(#1\right)}
  \newcommand{\mhom}[2]{\ensuremath{\mathrm{Hom}(#1,#2)}}
  \newcommand{\mshom}[3]{\ensuremath{\mathrm{Hom}_{#1}(#2,#3)}}
  \newcommand{\mend}[1]{\mathrm{End}(#1)}
  \newcommand{\msend}[2]{\ensuremath{\mathrm{End}_{#1}(#2)}}
  \newcommand{\mtr}[1]{\mathrm{Tr}\left(#1\right)}
  \newcommand{\mdet}[1]{\mathrm{det}\left(#1\right)}

  \newcommand{\ids}[0]{\mathrm{id}}
  \newcommand{\id}[1]{\ids_{#1}}

  \newcommand{\mdiag}[1]{\mathrm{diag}(#1)}
  \newcommand{\mexp}[1]{\mb{E}\left(#1\right)}
  \newcommand{\matr}[0]{M_d(\cc)}

  \newcommand{\spl}[0]{\langle}
  \newcommand{\spr}[0]{\rangle}
  \newcommand{\bva}[2]{e_{#1}^{\vphantom{#2}}}
  \newcommand{\bv}[1]{e_{#1}^{\vphantom{*}}}
  \newcommand{\dbv}[1]{e_{#1}^*}
  \newcommand{\vtn}[0]{V^{\otimes n}}
  \newcommand{\nullity}[1]{\mathrm{null}\left(#1\right)}
  \newcommand{\mtrsn}[2]{\mathrm{Tr}_{#1}\left(#2\right)}

  \newcommand{\mi}[1]{\bm{{#1}}}
  \newcommand{\indi}[0]{{\mi{i\vphantom{'}}}}
  \newcommand{\indj}[0]{{\mi{j\vphantom{'}}}}
  \newcommand{\indci}[0]{\mi{i'}}
  \newcommand{\indcj}[0]{\mi{j'}}
  
  \newcommand{\idns}[2]{\mathcal{I}_{#1}^{#2}}
  \newcommand{\idn}[0]{\idns{d}{n}}
  \newcommand{\parts}[2]{\mathcal{P}(#1,#2)}

  \newcommand{\grpu}[1]{\mathcal{U}(#1)}
  \newcommand{\grpgl}[1]{\mathrm{GL}(#1)}
  \newcommand{\grps}[1]{\mathcal{S}_{#1}}
  \newcommand{\charsn}[0]{\chi}
  \newcommand{\sncd}[0]{\cc_d[S_n]}
  \newcommand{\shortcentral}[1]{\mathrm{C}_{#1}}
  
  \newcommand{\grpalg}[1]{\cc\left[#1\right]}
  \newcommand{\gas}[1]{\grpalg{\grps{#1}}}
  \newcommand{\gau}[1]{\grpalg{\grpu{#1}}}
  \newcommand{\orb}[1]{\mathrm{orb}(#1)}
  \newcommand{\stab}[1]{\mathrm{stab}(#1)}
  \newcommand{\uint}[1]{\int_{\grpu{d}}#1 \, dU}  
  \newcommand{\uinti}[4]{\mathfrak{I}_d(#1,#2,#3,#4)}

  \newcommand{\ie}[0]{\textit{i.e.}}
  \newcommand{\eg}[0]{\textit{e.g.}}

  \newcommand{\fa}[0]{\ensuremath{\rightarrow} }

  \newcommand{\youngpos}[1]{\mathbb{Y}(#1)}
  
  \newcommand{\syt}[1]{\mathrm{SYT}(#1)}
  \newcommand{\cont}[1]{c(#1)}
  
  \newcommand{\charpi}[0]{\chi_\pi}
  \newcommand{\charpiarg}[1]{\charpi\left( #1 \right)}
  \newcommand{\charspecht}[0]{\chi_{\spv}}
  \newcommand{\charspechtarg}[1]{\charspecht\left( #1 \right)}
  \newcommand{\charweyl}[0]{\chi_{\wev}}
  \newcommand{\charweylarg}[1]{\charweyl\left( #1 \right)}
  \newcommand{\charjoint}[0]{\chi_{\pi\times\rho}}
  \newcommand{\charjointarg}[2]{\charjoint\left(#1,#2\right)}
  \newcommand{\charjointdecomp}[0]{\chi_{\spa \times\wea}}
  \newcommand{\charjointdecomparg}[2]{\charjointdecomp\left(#1,#2\right)}
  \newcommand{\spa}[0]{\pi_\lambda}
  \newcommand{\spv}[0]{S^\lambda}
  \newcommand{\spvm}[0]{S^\mu}
  \newcommand{\wea}[0]{\rho_\lambda}
  
  \newcommand{\wev}[0]{U^\lambda}
  \newcommand{\wevm}[0]{U^\mu}
  \newcommand{\integ}[0]{U(\indi,\indj,\indci,\indcj)}
  \newcommand{\triv}[1]{\trivnoarg(#1)}
  \newcommand{\trivnoarg}[0]{\mathcal{T}}
  
  \newcommand{\wein}[0]{\mathrm{Wg}}
  \newcommand{\weinarg}[1]{\mathrm{Wg}(#1)}

  \newcommand{\snmcparg}[1]{p^{#1}}
  \newcommand{\snmcp}[0]{\snmcparg{\lambda}}
  
  \newcommand{\cc}[0]{\mb{C}}
  \newcommand{\natnum}[0]{\mathbb{N}}
  \newcommand{\gagt}[1]{\ensuremath{\Phi\left(#1\right)}}

  \newcommand{\fullss}[0]{\substack{\lambda \vdash n \\ l(\lambda) \leq d}}
  \newcommand{\fulls}[0]{\lambda \vdash n,\, l(\lambda) \leq d}

  \newcommand\restr[2]{{
  \left.\kern-\nulldelimiterspace 
  #1 
  \vphantom{\big|} 
  \right|_{#2} 
  }}

\usepackage{csquotes} 
\usepackage[style=numeric,bibstyle=numeric,sorting=nty]{biblatex}
\usepackage[eulerchapternumbers=true]{classicthesis}

\clearscrheadfoot
\ohead[]{\headmark}
\rofoot[\pagemark]{\pagemark}
\lefoot[\pagemark]{\pagemark}

\usepackage{subfiles}

\title{Weingarten Calculus}
\author{Georg Köstenberger}

%% file: frontbackmatter/abstract.tex
\pdfbookmark[1]{Abstract}{Abstract}

\begingroup
\let\clearpage\relax
\let\cleardoublepage\relax
\let\cleardoublepage\relax

\chapter*{Abstract}
We consider the problem of computing the integral
\begin{equation}\label{eq-abs}
  \uint{u_{i_1j_1}\cdots u_{i_nj_n} \bar{u}_{i'_1j'_1} \cdots \bar{u}_{i'_{n'}j'_{n'}}},
\end{equation}
where the integration takes place with respect to the probability Haar measure on the unitary group, and the $u_{ij}$ denotes the $ij$-th entry of a unitary matrix $U$. We present a unified approach connecting classical results (\cite{hewittross}, \cite{weingarten}, \cite{samuel}), the explicit formula for the integral \eqref{eq-abs} given by B. Collins and P. Sniady \cite{collins2006} and subsequent works of various authors providing different points of view. Finally we are able to provide an explicit formula for the $2n$-th moment of the trace of a unitary Haar random matrix, generalizing a result of P. Diaconis \cite{diaconis}.
\vfill


\endgroup

\vfill

%% file: chapters/intro/intro.tex
\chapter{Introduction}\label{sec:1}
Let $\grpgl{d}$ denote the group of invertible, complex $d \times d$ matrices, and let $\grpu{d}$ be the subgroup of unitary matrices.
\graffito{For a proof of the existence and uniqueness of the Haar measure we refer to the paper G. K. Pedersen \cite{pedersen2004}.}
Recall that every locally compact topological group $G$ admits a regular Borel measure $\mu$, which is invariant under left translation \ie{}, $\mu(g\cdot X) = \mu(X)$ for all measurable sets $X$. This measure is unique up to a constant multiple. Any such measure is called \textit{left} Haar measure on $G$. Similarly there is  also a \textit{right} Haar measure on $G$. These two measure do not have to agree, but if they do $G$ is called \textit{unimodular}. 

The Haar measure induces a left-invariant integral on $G$,
\begin{equation*}
  \int_G f(hg) d\mu(g) = \int_G f(g) d\mu(g),
\end{equation*}
for all $h\in G$ and any Haar integrable function $f$ on $G$.

If $G$ is compact, then $\mu(G) < \infty$ and $G$ is unimodular. This is the reason why many arguments from the representation theory of finite groups carry over to compact groups almost verbatim. In many cases all we need to do is  exchanging $\frac{1}{|G|}\sum_{g\in G}$ with the normalized Haar measure on $G$ for the proofs to remain valid. 

As $\grpu{d}$ is compact, it admits a unique probability Haar measure.
The problem we are considering is the computation of the integrals with respect to this probability measure of the following kind

\begin{equation} \label{eq-basic-int}
\uint{u_{i_1j_1}\cdots u_{i_nj_n} \bar{u}_{i'_1j'_1} \cdots \bar{u}_{i'_{n'}j'_{n'}}},
\end{equation}

where the $u_{ij}$ denote the entries of a unitary matrix $U$, and $\bar{u}_{ij}$ denote the complex conjugate of said entries.

This at once allows us for example to define and compute an inner product on the algebra $\mathcal{A}$ of polynomial functions on $\grpu{d}$ \ie, the set of functions $f: \grpu{d} \fa \cc$ for which there exists a polynomial $p_f$ in $d^2$ variables with $f(U) = p_f(u_{11}, \dots, u_{dd})$.
As we will see, the integral  \eqref{eq-basic-int} is zero unless $n = n'$. This means that $\mathcal{A}$ admits a decomposition into homogeneous components \ie,

\begin{equation*}
  \mathcal{A} = \bigoplus_{k=0}^\infty \mathcal{A}^{(k)},
\end{equation*}

where $\mathcal{A}^{(k)}$ denotes the space of homogeneous polynomial functions of degree $k$ \cite{matnovak2009}.

Apart from answering natural questions about one of the classical groups, these integrals have many applications in random matrix theory and mathematical physics, especially quantum physics and quantum information theory. References for these applications can be found in \cite{linzhang}, \cite{matnovak2009}. Furthermore they have been used to derive formulas for the pseudoinverse of Gaussian matrices and the inverse of compound Wishart matrices \cite{colmatsaad}, and to describe integrals of Brownian motions on the classical groups \cite{intbrown}. 

The study of these integrals was initiated by physicists in the 1970s (\eg{} \cite{dewit77}). It was realized, that the integral can be expanded in terms of a sum over functions which only depend on the cycle structure of certain permutations. D. Weingarten was the first to study the asymptotic behavior of these functions \cite{weingarten}, hence these functions were later coined \textit{Weingarten functions} by B. Collins. Explicit formulas for the Weingarten functions were first derived by S. Samuel \cite{samuel} and B. Collins \cite{collins2002} under the assumption that $d \geq n$ and later by B. Collins and P. Sniady \cite{collins2006} in full generality. In the same paper B. Collins and P. Sniady also derived formulas for the integrals of polynomial functions on the orthogonal and symplectic groups. An alternative formula in terms of Jucys-Murphy elements was given by J. Novak \cite{novak2009} in 2009. This approach was generalized to the orthogonal group by P. Zinn-Justin \cite{zinnjustin} in 2010, and the computation of the Weingarten functions was seen to be equivalent to the computation of the Moore-Penrose inverse of a certain matrix.

There are several ways to tackle the problem of computing the integral \eqref{eq-basic-int}.
The first approach due to S. Samuel \cite{samuel} exploits certain symmetries of the integral, and leads to an \textit{ansatz} valid for $d\geq n$.
It is not obvious how to generalize this idea to a proof valid for all $n$ and $d$. 
On the other hand the method used by B. Collins and P. Sniady at first glance may look somewhat indirect, but leads to a valid general formula.
The methods used by B. Collins and P. Sniady, and thus all subsequent approaches as well, heavily depend on the Schur-Weyl duality and the double centralizer theorem. Furthermore we will make extensive use of both the classical as well as the modern approach to the representation theory of the symmetric group.

In order to motivate the general approach we compute the smallest possible example by hand in two different ways.

\begin{example}
  We seek to compute
\begin{equation*}
  \uint{|u_{11}|^2} = \uint{|u_{11}|^2}.
\end{equation*}
\end{example}

First we will use a direct approach. Let $U$ be a unitary Haar distributed $d\times d$ random matrix. All entries of $U$ have the same law. To see this note that permutation matrices are orthogonal and hence unitary. Left and right multiplication with permutation matrices can be used to exchange $u_{11}$ with any $u_{ij}$ we like. By the invariance of the Haar measure under such transformations all entries have the same distribution. Since the columns of $U$ are orthonormal we get
\begin{equation*}
  \uint{|u_{11}|^2} = \frac{1}{d}\sum_{i=1}^d \uint{|u_{1i}|^2} = \frac{1}{d}\uint{\sum_{i=1}^d |u_{1i}|^2} = \frac{1}{d}.
\end{equation*}

This method is direct, quick and only relies on elementary properties of the Haar-integral. However it is not obvious how to generalize this.

The following approach is in the spirit of B. Collins and P. Sniady. It is somewhat indirect, but easily lends itself to generalizations.

Let $D = \mathrm{diag}(1,0,\dots,0) = \bv{1}e_1^* \in \matr$ and note that
\begin{equation*}
  \begin{split}
    \uint{|u_{11}|^2} & = \uint{e_1^*U\bv{1} e_1^*U^*\bv{1}}\\
    & = \uint{\mtr{e_1^*U\bv{1} e_1^*U^*\bv{1}}} \\
    & = \uint{\mtr{\bv{1}e_1^*U\bv{1} e_1^*U^*}} \\
    & = \mtr{\uint{\bv{1}e_1^*U\bv{1} e_1^*U^*}}\\
    & = \mtr{D\uint{UDU^*}}.
  \end{split}
\end{equation*}
\graffito{A matrix that commutes with $\grpu{d}$ is a scalar multiply of the identity. This is a direct consequence of Schur's lemma (Lemma \ref{thm:schurs-lemma} and its corollary). Alternatively one can use the fact, that every matrix can be written as a linear combination of unitary matrices.}
By the invariance of the Haar measure we get that
\begin{equation*}
  W\left(\uint{UDU^*}\right)W^* = \uint{(WU)D(WU)^*} =  \uint{UDU^*},
\end{equation*}
for all $W \in \grpu{d}$.

Therefore $\uint{UDU^*}$ is a scalar multiple of the identity, say $\lambda \, \id{\cc^d}$. We get that
\begin{equation*}
  \begin{split}
    \lambda &= \frac{1}{d}\mtr{\uint{UDU^*}} \\
    &= \frac{1}{d}\uint{\mtr{UDU^*}} \\
    &= \frac{1}{d}\uint{\mtr{U^*UD}} \\
    & = \frac{1}{d}\uint{1} = \frac{1}{d}.
  \end{split}
\end{equation*}
Hence we compute
\begin{equation*}
  \uint{|u_{11}|^2} = \mtr{D\uint{UDU^*}} = \frac{1}{d}\mtr{D} = \frac{1}{d}.
\end{equation*}

The main benefit of the second approach is that it can be generalized quite easily.
Repeatedly applying the canonical isomorphism between tensor and Kronecker products 
\begin{equation*}
  A \otimes B \simeq \left(
    \begin{array}{c|c|c}
      a_{11}B & \,\cdots & a_{1d}B \\
      \hline
      \vdots &\,\ddots  & \vdots \\
      \hline
      a_{d1}B & \,\cdots & a_{dd}B \\
    \end{array} \right),
\end{equation*}
for $A, B \in \mend{V}$, to the matrix $U$, we see that the entries of $U^{\otimes n}$ and $(U^*)^{\otimes n}$ are $u_{i_1 j_1} \cdots u_{i_n j_n}$  and $\bar{u}_{i_1 j_1} \cdots \bar{u}_{i_n j_n}$ respectively.
Thus, if we try to select more than one pair of indices at once we end up with expressions of the form

\begin{equation*}
  \mtr{B \uint{U^{\otimes n} A (U^*)^{\otimes n}}},
\end{equation*}

for a suitable choice of $A$ and $B \in \matr^{\otimes n}$. Using this, the computation of the integrals reduces to the study of certain linear maps. In fact we will see, that the computation of the integral \eqref{eq-basic-int} is equivalent to the computation of an orthogonal projection $\mb{E}$ onto $\mend{\vtn}$. In order to describe them efficiently we introduce some notation.

\section{Notation}
Let $V = \cc^d$ and let $A$ be an element of $\mend{V^{\otimes n}} = \mend{V}^{\otimes n}$. We define
\begin{equation*}
  \mexp{A} = \uint{U^{\otimes n} A (U^*)^{\otimes n}},
\end{equation*}

where the integration is performed with respect to the unique probabilty Haar measure on $\grpu{d}$. 
As we will see in Proposition \ref{prop-e}, $\mb{E}$ is the orthogonal projection from $\mend{V}^{\otimes n}$ onto $\shortcentral{\grpu{d}}$, the centralizer of all $U^{\otimes n}$ for $U \in \grpu{d}$ in $\mend{V}^{\otimes n}$. The main takeaway from the above example is that the properties of $\mexp{A}$ and therefore the properties of the map $\mb{E}$ itself are central to our study. 

Additionally we introduce some notation to efficiently work with multiindices. 
Let $\idn = \set{(i_1, \dots,i_n) \mid 1 \leq i_k \leq d}$ be the set of $n$-tupels with entries in $1, \dots, d$. We denote multiindices with bold lowercase letters \eg{}, $\indi = (i_1, \dots,i_n)$, and note that the Kronecker delta for a pair of multiindices is given by
\begin{equation*}
  \delta_{\indi,\indj} = \prod_{k=1}^n \delta_{i_k,j_k} =
  \begin{cases}
    1 ~~~ \indi = \indj \\
    0 ~~~ \indi \neq \indj
  \end{cases}.
\end{equation*}
With this we write
\begin{equation*}
  \uinti{\indi}{\indj}{\indci}{\indcj} = \uint{u_{i_1j_1}\cdots u_{i_nj_n} \bar{u}_{i'_1j'_1} \cdots \bar{u}_{i'_{n'}j'_{n'}}},
\end{equation*}
for the integral.

This notation allows for a natural description of the integral in terms of tensor products.
Let $\set{\bv{i}}_{i=1}^d$ be a basis of $V$ and let $\set{\dbv{i}}_{i=1}^d$ be its dual basis.  We see that $\set{\bv{\indi}}_{\indi \in \idn}$ is a basis of $V^{\otimes n}$, where $\bv{\indi} = e_{i_1} \otimes \cdots \otimes e_{i_n}$. By setting
\begin{equation*}
  \bv{\indi}\dbv{\indj} = \bv{i_1}\dbv{j_1} \otimes \cdots \otimes \bv{i_n}\dbv{j_n}
\end{equation*}

we get that the $\bv{\indi}\dbv{\indj}$ for $\indi,\indj \in \idn$  form a basis of $\mend{V}^{\otimes n}$.
Therefore our problem is now to compute
\begin{equation*}
  \uinti{\indi}{\indj}{\indci}{\indcj} = \mtr{\bv{\indci}\dbv{\indi}\mexp{\bv{\indj}\dbv{\indcj}}} = \dbv{\indi}\mexp{\bv{\indj}\dbv{\indcj}}\bv{\indci}.
\end{equation*}

In other words, the entries of the matrix of $\mathbb{E}$ with respect to the standard basis of $\mend{\vtn}$ contain all the information we need. However the matrix of $\mb{E}$ contains a lot of redundant information as we will see, and computing $\mb{E}$ directly is quite a challenging task. Thus we do not explicitly compute $\mb{E}$. Instead we will use the representation theory of the symmetric group $\grps{n}$ and the Schur-Weyl duality for the unitary group $\grpu{d}$ to extract all the information we need from $\mb{E}$. In the next chapter we will develop the necessary theory for this.

%% file: chapters/intro/prelim.tex
\chapter{Preliminaries}
In this chapter we will briefly review basic notions from representation theory, the Haar measure and the Schur-Weyl duality.
We are mainly interested in the representation theory of finite and compact groups. Additionally we will need a few results from the representation theory of associative algebras. 
The representation theory of compact groups closely resembles the representation theory of finite groups. Moreover the algebras in question are group algebras of finite or compact groups and the representation theory of the group carries over to the group algebra. 
Thus we only consider the representation theory of finite groups.
This presentation is largely based on the exposition of the material given by J.-P. Serre \cite{serre96} with some results beeing taken from B. Simon \cite{bsim1996} and P. Etingof \cite{etingof2011}. For a treatment of the representation theory of compact groups and algebras, we refer to the books of D. Bump \cite{bump2013} and P. Etingof \cite{etingof2011} respectively.

\section{Representation theory}\label{section:reptheo}
Throughout this chapter we consider complex, finite dimensional representations, \ie, $V$ always denotes a complex, finite dimensional vector space.  

\begin{mdef}[Representation]
  Let $G$ be a finite (or compact) group and $\mathcal{A}$ an associative algebra. A representation of $G$ is a (continuous) homomorphism $\pi : G \fa \grpgl{V}$. A representation of $\mathcal{A}$ is a homomorphism $\pi : \mathcal{A} \fa \mend{V}$. 
  
  A subrepresentation of $\pi$ is a subspace of $V$ which is invariant under $\pi(x)$ for all $x$ in $G$ or $\mathcal{A}$. 

  If $V \neq 0$, $\pi$  is called irreducible, if its only subrepresentations are $0$ and $V$ itself. 
\end{mdef}

The group $G$ acts on a vector space $V$ via a representation $\rho$
\begin{equation*}
  g \cdot v := \rho(g)v.
\end{equation*}

We will sometimes call the vector space $V$ itself a representation of $G$ if there is no ambiguity regarding the associated action.  

Given two representations $\pi: G \fa \grpgl{V}$ and $\rho: G \fa \grpgl{W}$ we can construct a number of other representations.

They give rise to representations on $V \oplus W$ and $V \otimes W$, the latter via
\begin{equation*}
  g\cdot(v \otimes w) = (\pi(g)v) \otimes (\rho(g)w).
\end{equation*}

For every representation $\pi: G \fa \grpgl{V}$ we can define its \textit{dual representation} $\pi^* : G \fa \grpgl{V^*}$,  by $\pi^*(g) = \pi(g^{-1})^T$. This definition guarantees, that the dual paring $\spl\cdot ,\cdot \spr$ is invariant under the group action \ie,

\graffito{Since $\spl~ ,~ \spr$ is the dual pairing and not a complex scalar product we take the transpose and not the conjugate of $\pi(g^{-1})$.}
\begin{equation*}
  \langle \pi^*(g)v,\pi(g),w\rangle = \langle v, w\rangle.
\end{equation*}

We can also construct a representation of $G$ on $\mhom{V}{W}$, by setting
\begin{equation*}
  g\cdot T = \rho(g)\circ T \circ \pi(g^{-1}),
\end{equation*}
for $T \in \mhom{V}{W}$. This directly follows from unraveling the homomorphism $\mhom{V}{W} \simeq V^* \otimes W$ and the definition of the dual representation. 

\begin{mdef}[Intertwiner]
  An intertwiner between two representations $\pi:G \fa \mend{V}$ and $\rho:G \fa \mend{W}$ is a homomorphism $T: V \fa W$ which commutes with the action of $G$ \ie{}, $T(\pi(g)v) = \rho(g)T(v)$ or more concisely $T(g\cdot v) = g\cdot T(v)$.

  Two representations are said to be isomorphic, if there exists a bijective intertwiner between them.

  Furthermore we denote with \mshom{G}{V}{W} the space of all intertwining operators between $V$ and $W$. 
\end{mdef}

For every representation $U$ of $G$ we can define
\begin{equation*}
  U^G = \set{u \in U \mid g\cdot u = u,~ \forall g \in G},
\end{equation*}
the set of all elements of $U$ fixed under the action of $G$. In the special case of $U = \mhom{V}{W}$ we get $\mshom{G}{V}{W} = \mhom{V}{W}^G$.
\begin{lemma}[Schur's lemma]\label{thm:schurs-lemma}
  Let $V$ and $W$ be two representations of $G$ and $\phi \in \mshom{G}{V}{W}$. If $\phi \neq 0$ we have:
  \begin{enumerate}
  \item If $V$ is irreducible, $\phi$ is injective.
  \item If $W$ is irreducible, $\phi$ is surjective.
  \item If $V$ and $W$ are irreducible, $\phi$ is an isomorphism. 
  \end{enumerate}
\end{lemma}
\begin{proof}
  Note that $\mker{\phi}$ and $\mim{\phi}$ are subrepresentations of $V$ and $W$ respectively, and neither of them can be trivial since $\phi \neq 0$.
  
\end{proof}

\begin{corollary}
  Let $V$ be an irreducible representation of $G$ and let $\phi \in \msend{G}{V}$. Then $\phi = \lambda\,id_V$ for some $\lambda \in \cc$. 
\end{corollary}

\begin{proof}
  Since $\cc$ is algebraically closed the characteristic polynomial of $\phi$ has a root $\lambda$. $\phi - \lambda\,id_V$ is an intertwiner with nontrivial kernel and therefore by the previous lemma $\phi - \lambda\,id_V = 0$.
  
\end{proof}

\begin{corollary}\label{cor-delta}
  Let $U$ and $V$ be  irreducible representations of $G$. Then we have:
  \begin{enumerate}
  \item If $U \not \simeq V$, then $\mshom{G}{U}{V} = 0$.
  \item $\mdim{\msend{G}{V}} = 1$. 
  \end{enumerate}
\end{corollary}

\begin{lemma}\label{lem-subirrep}
  Every nonzero representation $\pi:G \fa \grpgl{V}$ has an irreducible subrepresentation. 
\end{lemma}
\begin{proof}
  Let $W = \set{\pi(g)v\mid g \in G}$ for some $v \neq 0$ in $V$. $W$ is a subrepresentation and $0 < dim\,W < dim\, V$. Either $W$ is already irreducible or it has a nontrivial subrepresentation $W'$, in which case we repeat the argument with $W'$ inplace of $W$.
  
\end{proof}

The standard theorems relating direct sums, tensor products and Hom-sets extend to representations. We have for example
\begin{equation*}
  \bigoplus_{i=1}^n \mshom{G}{V_i}{W} \simeq \mshom{G}{\bigoplus_{i=1}^n V_i}{W}.
\end{equation*}

In the case where $W$ is irreducible Schur's lemma implies that the dimension of $\mshom{G}{V}{W}$ is equal to the number of subrepresentations of $V$ isomorphic to $W$. 

\begin{mdef}[Group algebra]
  Let $G$ be a finite group.
  The group algebra $\cc[G]$ consists of all function $f:G \fa \mb{C}$ with the product given by convolution
  \begin{equation*}
    (fg)(x) = \sum_{u \in G} f(u)g(u^{-1}x)
  \end{equation*}
\end{mdef}

\begin{remark}
  $\cc[G]$ can also be thought of as the free vector space on $G$ over \mb{C}, with the product given by
  \begin{equation*}
    \left(\sum_{g \in G}\alpha_g g\right) \left(\sum_{h \in G}\beta_h h\right) = \sum_{g,h \in G}\alpha_g \beta_h gh.
  \end{equation*}

  A function $f: G \fa \mb{C}$ corresponds to the formal sum $\sum_{g \in G} f(g) g$.

  Note that $\cc[G]$ is itself a representation of $G$ with the action given by left multiplication, called the (left) regular representation. 
  We will use the functional notation $\delta_g$ whenever we want to emphasize that one should think of elements of $\cc[G]$ as function and we will use the free vector space notation when we are working with elements of $\cc[G]$ that are meant to act by left mulitplication. 
\end{remark}

\begin{remark}
  The definition of the group algebra can be adopted almost verbatim for compact groups. Let $dG$ denote the unique probability Haar measure on $G$ and define the convolution of $f,g: G \fa \cc$ to be
  \begin{equation*}
    (fg)(x) = \int_G f(u)g(x^{-1}u) dG(u).
  \end{equation*}
\end{remark}

\begin{theorem}[Maschke's theorem]
  Every subrepresentation $W$ of a  representation $V$ has a complement which is invariant under $G$, \ie{},  a subrepresentation. 
\end{theorem}

\begin{proof}
  There always exists an \textit{algebraic complement} $U'$ such that $V = W \oplus U'$, but $U'$ might not be invariant. Let $p$ be the projection along $U'$ onto $W$, and define
  \begin{equation*}
    P := \sum_{g \in G}\pi(g)p\pi(g^{-1}),
  \end{equation*}
  where $\pi$ is the action of $G$ on $V$. $P$ is a projection, $\mim{P} = W$ and hence $V = W \oplus \mker{P}$. Note that $\mker{P}$ is invariant and therefore $W$ has a complementary subrepresentation in $V$.
  
\end{proof}

\begin{remark}
  To adapt this proof to the compact case, we can use Weyl's unitary trick. Let $\pi: G \fa \grpgl{V}$ be a representation of a compact group $G$. If $V$ admits an inner product $\langle \cdot, \cdot \rangle$ invariant under the group action \ie{},
  \begin{equation*}
    \langle \pi(g) u, \pi(g) v \rangle = \langle v , u \rangle
  \end{equation*}
  for all $u, v \in V$, then the orthogonal complement with respect to this scalar product of any subrepresentation of $V$ is again a subrepresentation.
  Note that for an arbitrary scalar product $\langle \cdot, \cdot \rangle$ the scalar product
  \begin{equation*}
    (u,v) := \int_G \langle g \cdot u, g\cdot v \rangle dG
  \end{equation*}
  is well-defined by the compactness of $G$ and invariant.

\end{remark}

\begin{corollary}
  Every representation is a direct sum of irreducible representations. 
\end{corollary}

\begin{mdef}[Character]
  For a given representation $\pi : G \fa \grpgl{V}$ the function $\chi_\pi(g) := \mtr{\pi(g)}$ is called the character of $\pi$.
\end{mdef}

Characters of irreducible representations are called irreducible characters. We further note that $\chi_\pi \in \cc[G]$.

In the case of the dual representation we get
\begin{equation*}
  \chi_{\pi^*}(g) = \overline{\chi_\pi(g)}. 
\end{equation*}

Given \textit{some} action of $G$ on a finite set X, we can define a representation $\pi$ of $G$ by considering the free vector space $V = \set{e_x\mid x \in X}$ on $X$, setting $\pi(g)e_x = e_{g\cdot x}$. In this case the character of $\pi$ is equal to
\begin{equation*}
  \chi_\pi(g) = \mtr{\pi(g)} = \textit{fix}(g),
\end{equation*}
the number of elements of $X$ fixed under $g$.

Note that the characters of $V \oplus W$ and $V \otimes W$ are $\chi_V + \chi_W$ and $\chi_V\chi_W$ respectively.

We can define a scalar product on the group algebra $\cc[G]$ by setting
\begin{equation*}
  \spl \chi , \phi \spr  = \frac{1}{|G|} \sum_{g \in G}\overline{\chi(g)} \phi(g).
\end{equation*}

\begin{mprop}[Orthogonality of characters]
  Let $\pi : G \fa \grpgl{V}$ and $\rho: G \fa \grpgl{W}$ be two irreducible representations of $G$. Characters are orthogonal in the sense that
  \begin{equation*}
    \spl \chi_\pi,\chi_\rho \spr =
    \begin{cases}
      1 & V \simeq W \\
      0 & V \not\simeq W
    \end{cases}.
  \end{equation*}
\end{mprop}
\begin{proof}
  Note that for any representation $\sigma : G \fa \grpgl{U}$ the map
  \begin{equation*}
    \Phi = \frac{1}{|G|}\sum_{g\in G} \sigma(g) \in \msend{G}{U}
  \end{equation*}
  is a projection onto $U^G$.
  Furthermore we have
  \begin{equation*}
    \mdim{U^G} = \mtr{\Phi} = \frac{1}{|G|}\sum_{g \in G} \chi_\sigma(g).
  \end{equation*}

  If we set $U = \mhom{V}{W}$ and note that $\chi_{\mhom{V}{W}} = \overline{\chi_\pi(g)}\chi_\rho(g)$ we immediately get
  \begin{equation*}
    \mdim{\mshom{G}{V}{W}} = \frac{1}{|G|} \sum_{g \in G} \overline{\chi_\pi(g)} \chi_\rho(g) = \spl \chi_\pi , \chi_\rho \spr.
  \end{equation*}
  We conclude by applying Corollary \ref{cor-delta}. 
\end{proof}

Since characters are constant on conjugacy classes, we immediately get the following corollary.

\begin{corollary}
  The number of irreducible representations is less than or equal to the number of conjugacy classes. 
\end{corollary}

Hence we know that for any finite group the number of irreducible representation is itself finite. We denote with $\widehat{G}$ an indexing set of all irreducible representation of $G$. Any given representation $\pi: G \fa \grpgl{V}$ admits a decomposition
\begin{equation*}
  V \bigoplus_{\alpha \in \widehat{G}} n_a V_a
\end{equation*}
where $V_\alpha$ denotes an irreducible representation of $G$ and $n_\alpha$ denotes the number of times this irreducible representation appears in $V$ as a direct summand. The subresentation $n_\alpha V_\alpha$ is called the $\alpha$-isotypic component of $V$. 

Note that the character of the regular representation satisfies
\begin{equation*}
  \chi_{\cc[G]}(g) = |G|\delta_{g,e},
\end{equation*}
since multiplication in $\cc[G]$ has no fixed points with the exception of the neutral element.

Let $\widehat{G}$ be an indexing set for the irreducible representations of $G$ and let
\begin{equation*}
\cc[G] = \bigoplus_{\alpha \in \widehat{G}}n_\alpha V_\alpha
\end{equation*}
be the irreducible decomposition of $\cc[G]$. We have that
\begin{equation*}
  n_\alpha = \spl \chi_\alpha , \chi_{\cc[G]} \spr = \chi_\alpha(e) \frac{\chi_{\cc[G]}(e)}{|G|} = \mdim{V_\alpha}.
\end{equation*}

\begin{corollary}\label{ga-inj}
  The group algebra decomposes into
  \begin{equation*}
    \cc[G] \simeq \bigoplus_{\alpha \in \widehat{G}}\mdim{V_\alpha} \, V_\alpha \simeq \bigoplus_{\alpha \in \widehat{G}} \mend{V_\alpha}.
  \end{equation*}
\end{corollary}

\begin{proof}
  The second isomorphism warrants some further explanation, since it will become important later on.

  Consider the map $\phi: \cc[G] \fa \bigoplus_\alpha \mend{V_\alpha}$ defined by
  \begin{equation*}
    f \mapsto \left( \sum_{g \in G}f(g) \pi_\alpha(g) \right)_{\alpha \in \widehat{G}}
  \end{equation*}
  where $\pi_\alpha$ is the action on $V_\alpha$.

  $\phi(f) = 0$ implies that $\sum_g f(g) \pi_\alpha(g) = 0$ for all $\alpha \in \widehat{G}$. This in turn implies that $\sum f(g) \rho(g) = 0$ for any representation $\rho$. Applying this to the left regular representation yields $f = 0$. Hence this mapping is injective, and since the dimensions agree it is also bijective.
  
\end{proof}

Functions $f :G \fa \cc$ which are constant on conjugacy classes are called class functions. Note that characters are class functions.

For any representation $\rho: G \fa \grpgl{V}$ and any class function $f$ we define 
\begin{equation*}
  \rho_f = \sum_{g\in G} f(g)\rho(g) \in \msend{G}{V}.
\end{equation*}
If $V$ is irreducible the elements of \msend{G}{V} are just scalar multiples of the identity map and we get that
\begin{equation*}
  \rho_f = \frac{\mtr{\rho_f}}{\mdim{V}}~\id{V} = \frac{|G|}{\mdim{V}}\spl \chi_\rho^*, f \spr~ \id{V}.
\end{equation*}

\begin{mprop}
  The irreducible characters are an orthonormal basis of the class functions. 
\end{mprop}
\begin{proof}
  We show that $\spl \chi_\alpha^*\spr^{\bot} = 0$. For $f \in \spl \chi_\alpha^*\spr^{\bot}$ we have $\rho_f = 0$ for any irreducible representation $\rho$ and hence for any representation. Applying this to the left regular representation yields $f = 0$.
  
\end{proof}

\begin{corollary}\label{thm:conj-class}
  The number of irreducible characters is equal to the number of conjugacy classes. 
\end{corollary}

\begin{mdef}[minimal central projection]
  A projection (\ie{} an idempotent) $p$ in $A$ is called minimal if for every pair of projections $r$ and $q$ we have that $p = r + q$ implies $p = 0$ or $q = 0$. The projection $p$ is called minimal central if this is true for all projections $q, r$ in the center of $A$.  
\end{mdef}

\graffito{For a more detailed treatment of minimal projections we  refer to the book by B. Simon \cite{bsim1996}.}
Let $\rho: G \fa \grpgl{W}$ be an arbitrary representation of $G$ and let $\pi: G \fa \grpgl{V}$ be an irreducible one. Then
\begin{equation*}
  \frac{\mdim{V}}{|G|}\rho_{\chi_\pi^*}
\end{equation*}
is the minimal central projection onto the $\pi$-isotypic component of $\rho$. 

\begin{mdef}[Partition]
  A partition of an integer $n\in \natnum$ is a tuple $\lambda = (\lambda_1,\dots,\lambda_l)$ such that $\lambda_1 \geq \cdots \geq \lambda_l > 0$ and $\lambda_1 + \cdots + \lambda_l=n$. The number $l$ is called the length of $\lambda$ and denoted by $l(\lambda)$. If $\lambda$ is a partition of $n$ we write $\lambda \vdash n$. 
\end{mdef}

Note that there is a bijection between partitions of $n$ and the conjugacy classes of $\grps{n}$.
This is due to the following considerations.
Two permutations are conjugate if and only if they have the same cycle structure.
Furthermore disjoint cycles commute in $\grps{n}$. This allows us to decompose any permutation $\sigma$ into a composition of disjoint cycles of descending length $\sigma_1\cdots\sigma_l$.
The lengths of the cycles $l(\sigma_i)$ sums up to $n$.
This yields a partition $\lambda(\sigma) = (l(\sigma_1), \dots, l(\sigma_l))$ of $n$. 

To each partition $\lambda$ we can associate its Young diagram. A Young diagram is a finite collection of boxes arranged in a top-left-justified rectangle, such that the $i$-th row  contains $\lambda_i$ boxes. The box in the $i$-th row and $j$-th column will be denoted by the pair $(i,j)$. We define $\youngpos{n}$ to be the set of all Young diagrams associated to partitions of $n$ and $T: \grps{n} \fa \youngpos{n}$ the map mapping a permutation to the Young diagram of its conjugacy class. For later use we define for each $\mu \in \youngpos{n}$ the element $C_\mu \in \gas{n}$ to be
\begin{equation*}
  C_\mu = \sum_{T(\sigma) = \mu} \delta_\sigma,
\end{equation*}
\ie, the characteristic function of the conjugacy class of $\mu$. As for partitions let $l(\mu)$ denote the number of rows of $\mu \in \youngpos{n}$.  

In view of Corollary \ref{thm:conj-class} this means that the irreducible representations of $\grps{n}$ are indexed by partitions.
They are called Specht modules and will be denoted by $\spv$.
The associated action of $\grps{n}$ on $\spv$ will be denoted by $\spa$ and we will write $\charspecht$ for its character.

\marginpar{%
  \includegraphics[width=\marginparwidth]{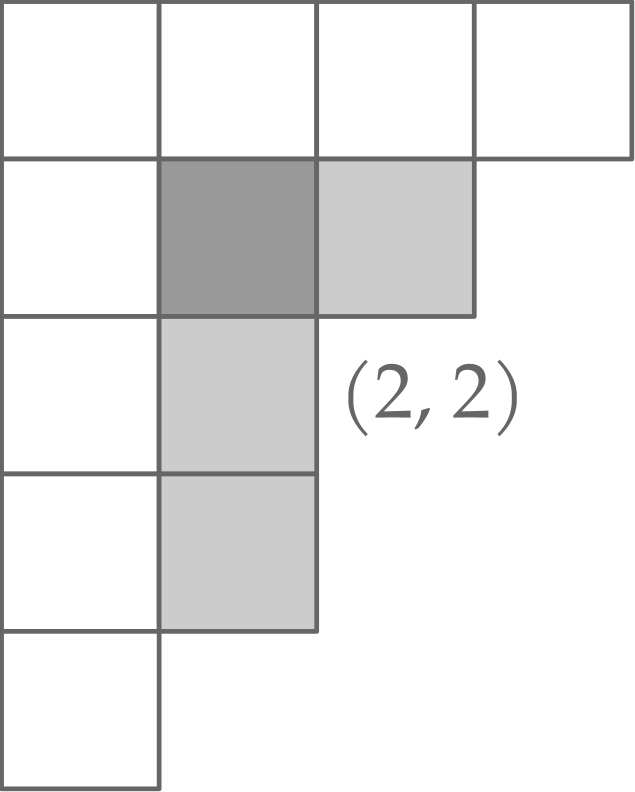}
  \captionof*{figure}{Young diagram of the partition $\lambda=(4,3,2,2,1)$ with the hook of the box $(2,2)$ in grey.}
}

In order to compute the dimensions of the representations that are going to appear, we need one more concept related to Young diagrams, namely the \textit{hook} of a box $(i,j)$.
For a given Young diagram $\mu$ the hook $H_\mu(i,j)$ consists of all boxes below or to the right of the box $(i,j)$.
The \textit{hook length} $h_\mu(i,j)$ is the number of boxes in $H_\mu(i,j)$ (the box $(i,j)$ also belongs to the hook). 

The classical approach to the representation theory of $\grps{n}$ uses Young Tableaux/Young Symmetrizers, as detailed in W. Fulton and J. Harris \cite{fh2004}.
Additionally there is an alternative approach due to A. Okounkov and A. Vershik \cite{rep-sym-new}.
Discussing the details of this alternative approach is beyond the scope of this work, however we will briefly touch upon a certain construction used by Okounkov and Vershik. The Jucys-Murphy elements $J_i$ are special elements of the group algebra $\gas{n}$, independently introduced by G. E. Murphy \cite{murphy} and A.-A. A. Jucys \cite{jucys}. They are defined as sums of the Coxeter generators of $\grps{n}$: 
\begin{equation*}
  \begin{split}
    J_1 &= 0 \\
    J_2 &= (1,2) \\
    J_3 &= (1,3) + (2,3) \\
    \vdots& \\
    J_n &= (1,n) + (2,n) + \cdots + (n-1,n).
  \end{split}
\end{equation*}
Let $\gas{k-1} \subseteq \gas{n}$ be the subalgebra consisting of all permutations of $\set{1,\dots,k-1}$. Note that for any $\sigma \in \gas{k-1}$ we have
\begin{equation*}
  \sigma J_k \sigma^{-1} = \sum_{i=1}^{k-1}\sigma (i,k) \sigma^{-1} = \sum_{i=1}^{k-1}(\sigma(i),k) = J_k. 
\end{equation*}
Thus any two Jucys-Murphy elements commute with each other and $J_k$ commutes with the subalgebra $\gas{k-1}$ of $\gas{n}$.
In fact the subalgebra $\cc[J_2,\dots,J_n] \subseteq \gas{n}$ is a maximal commutative subalgebra, known as \textit{Gelfand-Tsetlin} algebra. 
The Jucys-Murphy elements possess a number of very useful properties.
However we will only need the following result due to Jucys and a rough estimate of their \textit{eigenvalues}.

\begin{theorem}[Jucys, \cite{jucys}]\label{thm:jucys}
  For $r \in \{1, \dots, n\}$, let
  \begin{equation*}
    e_r(J_1, \dots, J_n) = \sum_{1 \leq i_1 < \cdots < i_r \leq n} J_{i_1} \cdots J_{i_r}
  \end{equation*}
  be the elementary symmetric polynomials evaluated in the Jucys-Murphy elements. Then
  \begin{equation*}
    e_r(J_1, \dots, J_n) = \sum_{\substack{\mu \in \youngpos{n} \\ l(\mu) = n-r}} C_\mu.
  \end{equation*}
\end{theorem}

Note that this means that $e_r(J_1, \dots, J_n)$ is the characteristic function of all permutations with exactly $n-r$ cycles.

We can think of every element in $\gas{n}$ as a vector space homomorphism by considering the left regular representation (\ie{} acting by left multiplication) of $\gas{n}$ on itself. Note that $\gas{n}$ has a basis indexed by permutations. Hence the action of $\sigma \in \gas{n}$ corresponds to a $n!\times n!$ matrix $L(\sigma)$  with complex entries. The multiplication in $\gas{n}$ mimics the multiplication of these matrices. When we speak of eigenvalues of $\sigma \in \gas{n}$ we mean the eigenvalues of the matrix corresponding to left multiplication by $\sigma$.

Note that the matrix $L(J_k)$ of any Jucys-Murphy element $J_k$ is symmetric, since all the matrices corresponding to $(i,k)$ are. Furthermore the entries of $L(J_k)$ are either $0$ or $1$ and all columns (and hence all rows) of $L(J_k)$ sum to $k$. The next lemma proves that all eigenvalues of $J_{k+1}$ lie in the interval $[-k,k]$.

\graffito{The lemma and its application to the eigenvalues of the Jucys-Murphy elements is due to an answer of D. E. Speyer on the \texttt{mathoverflow} stackexchange \cite{speyer-jme-overflow}.}

\begin{lemma}\label{lemma:jucys-ev}
  Let $A$ be a symmetric matrix with nonnegative entries, whose rows and columns all sum to $k$. Then the eigenvalues of $A$ lie in the interval $[-k,k]$. 
\end{lemma}

\begin{proof}
  For any vector $v$ we have
  \begin{equation*}
    \spl v, Av \spr = k \sum_i v_i^2 - \sum_{i<j}A_{ij}(v_i-v_j)^2 \leq k \sum_i v_i^2 = k\spl v, v\spr.
  \end{equation*}
  Hence the Reyleigh quotient is atmost $k$, which implies the statement. 
\end{proof}

In order to discuss P. Zinn-Justin's approach, we need one more classical concept from the representation theory of $\grps{n}$, namely the notion of a Standard Young Tableau. Let $\lambda \vdash n$ be a partition of $n$, and $T$ the Young diagram corresponding to $\lambda$. We can fill the boxes of $T$ with the numbers $1, \dots, n$ to obtain a \textit{Young tableau}. If the entries in each row and column are increasing, then the Young tableau is called a \textit{standard Young tableau}. We write $\syt{\lambda}$ for the set of all standard Young tableaux of shape $\lambda$. This allows us to introduce the set of Young's orthogonal idempotents $e_T$ for a standard Young tableau $T$ with $n$ boxes. They are completely characterized by the following properties 
\begin{equation}
  \label{eq:young-idem-1}
  e_Te_S = \delta_{TS}e_T ~~~ \text{and} ~~~ \sum_{T:|T| = n} e_T = 1.
\end{equation}
Furthermore they diagonalize the Jucys-Murphy elements,
\begin{equation}
  \label{eq:young-idem-2}
  J_ke_T = e_TJ_k = \cont{T_k}e_T, ~~~ \text{for } k = 1, \dots, n,
\end{equation}
where $\cont{T_k}$ is the \textit{content} of the box of $T$ labelled $k$, defined by $\cont{T_k} = j-i$ if the box $(i,j)$ has label $k$. 

Finally we can define the central idempotent $\snmcp$ associated to $\lambda \vdash n$,
\begin{equation}
  \label{eq:young-central-idem}
  \snmcp  = \sum_{T\in \syt{\lambda}} e_T.
\end{equation}
From the definition of $e_T$ it is immediate that
\begin{equation*}
 \snmcp \snmcparg{\mu} = 0 ~~~ \text{for} ~~~ \mu \neq \lambda, ~~~ \text{and} ~~~ (\snmcp)^2 = \snmcp.
\end{equation*}

One can show, that $\snmcp$ is nothing else but the minimal central projection onto the $\lambda$-isotypic component of $\pi$. Hence it can be written as
\begin{equation*}
  \snmcp = \frac{\mdim{\spv}}{n!} \charspecht.
\end{equation*}

The irreducible representations $\wea: \grpu{d} \fa \grpgl{\wev}$ of the unitary group $\grpu{d}$ can be indexed by partitions $\lambda \vdash n$, such that $l(\lambda) \leq d$. The $\wev$'s are called \textit{Weyl modules} \cite{fh2004}.

The next theorem will allow us to compute the Weingarten functions in practice.

\begin{theorem}[Hook length formula] \label{thm:hook-length}
  The dimensions of the irreducible representations of $\grps{n}$ and $\grpu{d}$ are given by
  \begin{equation*}
    \begin{split}
      \mdim{\spv} & = \frac{n!}{\prod_{(i,j)\in \lambda}h_\lambda(i,j)},\\
      \mdim{\wev} & = \prod_{(i,j) \in \lambda} \frac{d+j-i}{h_\lambda(i,j)} = \prod_{1\leq i < j \leq d} \frac{\lambda_i - \lambda_j + j -i}{j - i}.
    \end{split}
  \end{equation*}
\end{theorem}

Next we will introduce two group actions on $\idn$, one for $\grps{n}$ and one for $\grps{d}$. These actions will lift to representations of $\grps{n}$ and $\grps{d}$ on $\vtn$. Furthermore we will introduce a representation of $\grpu{d}$ on $\vtn$. The relationship between the representations of $\grps{n}$ and $\grpu{d}$ will be one of the cornerstones of our theory. 

The action of $\grps{n}$ is denoted with $\sigma(\indi)$ and defined by \textit{permuting} the entries \ie, $\sigma(\indi) = (i_{\sigma^{-1}(1)},\dots,i_{\sigma^{-1}(n)})$. This lifts to a representation $\pi:\grps{n} \fa \vtn$ of $\grps{n}$ via $\pi(\sigma)\bv{\indi} = \bv{\sigma(\indi)}$. Thus for arbitrary $v_1, \dots ,v_n \in V$ we have
\begin{equation*}
  \pi(\sigma): v_1 \otimes \cdots \otimes v_n \mapsto v_{\sigma^{-1}(1)} \otimes \cdots \otimes v_{\sigma^{-1}(n)}.
\end{equation*}

\begin{remark}
  We have to make sure that $(\sigma\tau)(\indi) = \sigma(\tau(\indi))$ in order for it to be a group action. If we were to define the action of $\sigma \in \grps{n}$ for an arbitrary $\indi \in \idn$ by $\sigma(\indi) = (i_{\sigma(1)},\dots,i_{\sigma(n)})$ we would end up with the right action
  \begin{equation*}
    \begin{split}
      \sigma(\tau(\indi)) &= \sigma(i_{\tau(1)}, \dots , \i_{\tau(n)}) \\
      & = \sigma(j_1, \dots, j_n) \\
      & = (j_{\sigma(1)},\dots, j_{\sigma(n)}) \\
      & = (i_{\tau(\sigma(1))},  \dots, i_{\tau(\sigma(n))}) \\
      & = (\tau\sigma)(\indi),
    \end{split}
  \end{equation*}
  where $j_i = i_{\tau(i)}$. Permuting according to the inverse of $\sigma$ yields a left action. 
  
\end{remark}

\begin{remark}
  In view of the isomorphism $\vtn \simeq \cc^n \otimes V$ we see that $\pi$ acts on the $\cc^n$ component, while $\rho$ acts on $V$,
  \begin{equation*}
    \begin{split}
      \pi(\sigma)& = \sigma \otimes \id{V}, \\
      \rho(U) &= \id{\cc^n} \otimes U.
    \end{split}
  \end{equation*}
\end{remark}


The action of $\grps{d}$ is denoted with $\rho[\indi]$ and \textit{exchanges} the entries of multiindices \ie{}, $\rho[\indi] = (\rho(i_1), \dots, \rho(i_n))$. This again yields an action of $\grps{d}$ on $\vtn$, that is best understood in a slightly different fashion than $\pi$. For any $\rho \in \grps{d}$ define $S_\rho$ to be the linear map
\begin{equation*}
  S_\rho  = \sum_{i=1}^d \bv{\rho(i)}\dbv{i}.
\end{equation*}

Clearly we have $S^{\otimes n}_\rho\bv{\indi} = \bv{\rho[\indi]}$. Furthermore since $S_\rho$ is unitary, the actions of $\grps{n}$ and $\grps{d}$ commute. Note that $S_\rho^{-1} = S_\rho^T = S_\rho^* = S_{\rho^{-1}}$.

The representation $\rho: \grpu{d} \fa \vtn$ of the unitary group is given by the tensor power \ie{}, $\rho(U) = U^{\otimes n}$. Thus for any $v_1, \dots, v_n \in V$ we have  
\begin{equation*}
  \rho(U): v_1 \otimes \cdots \otimes v_n \mapsto (Uv_1) \otimes \cdots \otimes (Uv_n).
\end{equation*}

Note that $\pi(\sigma)$ commutes with every $U^{\otimes n}$ for $U \in \grpu{d}$. In fact $\pi(\sigma)$ commutes with the diagonal action of $\mend{V}$ on $\mend{V}^{\otimes n}$ \ie{}, with every $M^{\otimes n}$ for $M \in \mend{V}$. Furthermore note that $\pi(\sigma)^* = \pi(\sigma)^T = \pi(\sigma^{-1}) = \pi(\sigma)^{-1}$. Let $\charsn(\sigma) = \mtr{\pi(\sigma)}$ be the character of $\pi$. This map will be of central importance, thus we compute it explicitly.
\begin{lemma}\label{lemma:charpi}
  The character of $\pi$ is given by
  \begin{equation*}
    \charsn(\sigma) = d^{\#\sigma},
  \end{equation*}
  where $d = \mdim{V}$ and $\#\sigma$ denotes the number of disjoint cycles of $\sigma$ in the canonical factorization. 
\end{lemma}
\begin{proof}
  As $\pi(\sigma)$ is a permutation matrix, its trace is given by the number of fixed points of the group action. In order for a multiindex $\indi$ to be fixed under this action, the elements of every cycle of $\sigma$ have to be identical. There are $d$ elements with which we can fill any cycle, and $\#\sigma$ cycles in total, yielding exactly $d^{\#\sigma}$ elements fixed under the action of  $\sigma$.   
\end{proof}
Since the actions $\pi$ and $\rho$ commute they give rise to a joint action of $S_n \times U(d)$ on $\vtn$, which is the subject of the next section. 

\section{Schur-Weyl duality}\label{sec:schur-weyl}
Let $\mathcal{A} = \spl \pi(\sigma) \mid \sigma \in \grps{n}\spr$ and $\mathcal{B} = \spl U^{\otimes n} \mid U \in \grpu{d}\spr$ be the subalgebras of $\mend{V^{\otimes n}}$ generated by the actions of $\grps{n}$ and $\grpu{d}$ respectively. Note that $\mathcal{A} = \pi(\gas{n})$ and $\mathcal{B} = \rho(\gau{d})$. As it will turn out, the centralizer $\shortcentral{\grpu{d}} = \msend{\mathcal{B}}{V^{\otimes n}}$ of $\mathcal{B}$ and the question of how exactly we can decompose $\mathcal{A}$ into images of irreducible representations of $\grps{n}$ will be of central importance. These questions are the subject of the double centralizer theorem.   
The following version is taken from Etingof et al. \cite{etingof2011}. 
\begin{theorem}[Double Centralizer Theorem]\label{thm:double-cent}
  Let $A$ and $B$ be two subalgebras of $\mend{V}$. If $A$ is semisimple and $B=\msend{A}{V}$ then:
  \begin{enumerate}
  \item $B$ is semisimple.
  \item $A = \msend{B}{V}$.
  \item As a representation of $A \otimes B$, $V$ decomposes into
    \begin{equation}
      \label{eq:dc:decomp}
      V = \bigoplus_{i \in I} V_i \otimes W_i,
    \end{equation}
    where $U_i$ are all the irreducible representations of $A$ and $W_i$ are all the irreducible representations of $B$. 
  \end{enumerate}
\end{theorem}
\begin{proof}[Proof (Sketch)]
  Since $A$ is semisimple it has an isotypical decomposition, $A = \bigoplus_{i\in I} \mend{V_i}$. Furthermore $V \simeq \bigoplus_{i\in I} V_i  \otimes \mshom{A}{V_i}{V}$, since the characters of the two representations agree. By Schurs's lemma $B = \msend{A}{V} \simeq \bigoplus_{i\in I} \mend{W_i}$ with $W_i = \mshom{A}{V_i}{V}$. This implies all statements of the theorem. 
\end{proof}

Applying this to our objects of interest, $\mathcal{A} = \sigma(\gas{n})$ and $\mathcal{B} = \rho(\gau{d})$, we get the following classical result. The proof is taken from R. Goodman and N. R. Wallach \cite{goodmanwallach}.

\begin{theorem}\label{thm:cent}
  $\mathcal{A}$ and $\mathcal{B}$ are centralizers of each other \ie{},
  \begin{equation}
    \label{eq:sw-central}
    \msend{\mathcal{A}}{V^{\otimes n}} = \mathcal{B} ~~~ \text{and} ~~~\msend{\mathcal{B}}{V^{\otimes n}} = \mathcal{A}.
  \end{equation}

\end{theorem}
\begin{proof}
  Let $B = (b_{\indi \indj}) \in \mend{V^{\otimes n}}$. Since
  \begin{equation*}
    B\pi(\sigma)e_\indj = Be_{\sigma(\indj)} = \sum_{\indi} b_{\indi \sigma(\indj)}e_\indi
  \end{equation*}
  and
    \begin{equation*}
    \pi(\sigma^{-1})Be_\indj  = \sum_{\indi} b_{\indi \indj}e_{\sigma^{-1}(\indi)}= \sum_{\indi} b_{\sigma(\indi) \indj}e_{\indi},
  \end{equation*}
  we have $B \in \msend{\mathcal{A}}{V^{\otimes n}}$ if and only if $b_{\indi,\indj} = b_{\sigma(\indi) \sigma(\indj)}$ for all $\indi,\indj \in \idn$ and $\sigma \in \grps{n}$.

  Consider the nondegenerate bilinear form $\spl X, Y \spr := \mtr{XY}$. We will show that the restriction of this form to $\msend{\mathcal{A}}{V^{\otimes n}}$ is still non-degenerate. To this end consider the projection $P$ from $\mend{V^{\otimes n}}$ onto $\msend{\mathcal{A}}{\vtn}$ explicitly given by
  \begin{equation*}
    P(X) = \frac{1}{n!}\sum_{\sigma \in \grps{n}} \pi(\sigma) X \pi(\sigma^{-1}).
  \end{equation*}
  For any $B \in \msend{\mathcal{A}}{\vtn}$ and $X \in \mend{\vtn}$ we get
  \begin{equation*}
    \spl P(X), B \spr = \frac{1}{n!} \sum_{\sigma \in \grps{n}} \mtr{\pi(\sigma) X \pi(\sigma^{-1}) B} = \spl X , B \spr. 
  \end{equation*}
  Thus $\spl A, B \spr = 0$  for all $A \in \msend{\mathcal{A}}{\vtn}$ implies $\spl X , B \spr = 0$ for all $X \in \mend{\vtn}$ and by the non-degeneracy of $\spl\cdot,\cdot\spr$ on $\mend{\vtn}$ we get that $B = 0$. Therefore it is sufficient to show that $\msend{\mathcal{A}}{\vtn} \cap \mathcal{B}^{\perp} = 0$ in order to prove that $\msend{\mathcal{A}}{\vtn} = \mathcal{B}$.

  In fact there is room for a slight generalization. We do not need to restrict ourselves to unitary transformations. The following argument applies to both, $\grpgl{d}$ and $\grpu{d}$.

  Note that for any $G \in \grpgl{d}$ we have $G^{\otimes n} = (g_{\indi\indj}) = (g_{i_1 j_1} \cdots g_{i_n j_n})$. For a fixed $B \in \msend{\mathcal{A}}{\vtn} \cap \mathcal{B}^{\perp}$ this yields
  \begin{equation*}
    \spl B, G \spr = \sum_\indi e_\indi^* B G \bva{\indi}{T} = \sum_{\indi,\indj} b_{\indj \indi} g_{\indj \indi}.
  \end{equation*}
  Associate to $B$ a polynomial function $f_B$ on $\mend{V}$,
  \begin{equation}
    \label{eq:gl-poly}
    f_B(X) = \sum_{\indi, \indj} b_{\indi\indj}x_{\indi\indj} = \sum_{\indi, \indj} b_{\indi\indj} x_{i_1j_1} \cdots x_{i_n j_n},
  \end{equation}
  for $X \in \mend{V}$. Note that $f_B$ vanishes on  $\grpgl{d}$ and since $\grpgl{d}$ is dense in $\mend{V}$ we have $f_B \equiv 0$. Plugging this into equation \eqref{eq:gl-poly} we get
  \begin{equation}
    \label{eq:gl-poly-0}
    \sum_{\indi, \indj} b_{\indi \indj} x_{i_1j_1} \cdots x_{i_n j_n} = 0.
  \end{equation}
  Let $\indi$ and $\indj$ be fixed and define $X \in \mend{V}$ according to
  \begin{equation*}
    x_{ij} =
    \begin{cases}
      1 & \text{if } i \in \set{\indi} \text{ and } j\in \set{\indj}\\
      0 & \text{otherwise}\\
    \end{cases}
  \end{equation*}
  Note that $x_{\indci \indcj} = 1$ if and only if $\indci = \sigma(\indi)$ and $\indcj = \sigma(\indcj)$ for some $\sigma \in \grps{n}$. Since $b_{\sigma(\indi)\sigma(\indj)} = b_{\indi,\indj}$ for all $\sigma \in \grps{n}$ this implies that $B = 0$. 
\end{proof}

\begin{corollary}
  With $\mathcal{A}$ and $\mathcal{B}$ as before we have
  \begin{equation*}
    \msend{\mathcal{A}}{\vtn} = \spl G^{\otimes n} \mid G \in \grpgl{d}\spr = \spl T^{\otimes n} \mid T \in \mend{V}\spr =  \mathcal{B}
  \end{equation*}
\end{corollary}

Combining Theorem \ref{thm:double-cent} and \ref{thm:cent} yields the main part of the Schur-Weyl duality. We do not prove that the irreducible representations of $\grpu{d}$ are indexed by partitions $\lambda \vdash n$ with length $l(\lambda) \leq d$, as this would require the notion of a weight space and the so called \textit{theorem of the highest weight}. For an exposition of this fact, we refer to \cite{fh2004}.

\begin{theorem}[Schur-Weyl Duality]\label{thm:swd}
  The action of $\grps{n}\times \grpu{d}$ is multiplicity free, \ie{}, no irreducible representation of $\grps{n} \times \grpu{d}$ appears twice. The decomposition into irreducible components is given by
   
  \begin{equation}
    \label{eq:sw-decomp}
    \vtn = \bigoplus_{\fullss} \spv \otimes \wev.
  \end{equation}
\end{theorem}




%% file: chapters/main/main.tex
\chapter{Integration on unitary groups}
In this chapter we are going to discuss several ways to compute the integrals $\uinti{\indi}{\indj}{\indci}{\indcj}$ and related quantities. Furthermore we are going to derive closed form expressions for special choices of $\indi$, $\indj$, $\indci$ and $\indcj$. We start with some elementary properties of the integral $\uinti{\indi}{\indj}{\indci}{\indcj}$, that do not require any sophisticated theory. Then we will discuss the general approach to the computation of $\uinti{\indi}{\indj}{\indci}{\indcj}$ due to B. Collins and P. Sniady and its consequences. Finally we will consider alternative approaches using the theory of Jucy-Murphy elements and Moore-Penrose inverses. 

\section{Elementary properties}
\begin{mprop}\label{prop:int-vanish}
  The integral vanishes unless $n = n'$. 
\end{mprop}
\begin{proof}
  Let $k$ be an integer such that $k$ does not divide $n-n'$ and let $u$ be a $k$-th root of unity. Note that $u \, \id{V}  \in \grpu{d}$. Since the Haar measure is invariant under the group action we get
  \begin{equation*}
    u^{n-n'}\int_{U(d)}\integ dU = \int_{U(d)}\integ dU,
  \end{equation*}
  but $u^{n-n'} \neq 1$ by our choice of $k$ and hence the integral vanishes.
  
\end{proof}

Variations of this argument will allow us to prove many properties of the integral, but first we consider the basic properties of the map $\mb{E}$.

Since the expectation of a random matrix $U$ is defined componentwise ($\mexp{U}$ corresponds to the matrix $(\mexp{U_{ij}})_{ij}$) the previous theorem already has some interesting consequences. 

\begin{mprop}
  For $k \geq 1$ we have
  \begin{enumerate}
  \item $\mb{E}(U^k) = 0$.
  \item $\mb{E}(\mtr{U}^k) = \mb{E}(\mtr{U^k}) = 0$.
  \item $\mb{E}(\mdet{U}^k) = 0$. 
  \end{enumerate}
\end{mprop}

The next proposition is concerned with the basic properties of the map $\mexp{A} = \uint{U^{\otimes n} A (U^*)^{\otimes n}}$. 
\begin{mprop}\label{prop-e}
  The map $\mb{E}$ is an orthogonal projection onto $\shortcentral{\grpu{d}}$ with respect to the scalar product $\spl A, B \spr = \mtr{A^*B}$ having the following properties:
  \begin{enumerate}
  \item $W^{\otimes n}\mexp{A}(W^*)^{\otimes n} = \mexp{A} = \mexp{W^{\otimes n}A(W^*)^{\otimes n}}$ for $W\in \grpu{d}$.
  \item $\mtr{A} = \mtr{\mexp{A}}$.
  \item $\mexp{XAY} = X \mexp{A} Y$  for all $X,Y \in \shortcentral{\grpu{d}}$.
  \item $\mtr{A\mexp{B}} = \mtr{\mexp{A}B}$.
  \item $\mexp{\bv{\indi}\dbv{\indj}} = 0$ unless $\indj = \sigma(\indi)$ for some $\sigma \in \grps{n}$.
  \item $\mexp{\bv{\indi}\dbv{\indj}} = \mexp{\bv{\rho[\indi]}\dbv{\rho[\indj]}}$ for all $\rho \in \grps{d}$. 
  \end{enumerate}
\end{mprop}
\begin{proof}
  We start by noting, that $\mexp{A} = A$ for every $A$ in the centralizer of the maps $U^{\otimes n}$ for $U \in \grpu{d}$, and hence $\mim{\mathbb{E}} \subseteq \shortcentral{\grpu{d}}$. By the invariance of the Haar measure we get
  \begin{equation*}
    \begin{split}
      W^{\otimes n}\mexp{A}(W^*)^{\otimes n} &= \uint{(WU)^{\otimes n} A ((WU)^*)^{\otimes n}} \\
      & = \uint{U^{\otimes n}A (U^*)^{\otimes n}} = \mexp{A}\\
      & = \uint{(UW)^{\otimes n}A(UW)^{\otimes n}} = \mexp{W^{\otimes n}A(W^*)^{\otimes n}}.
    \end{split}
  \end{equation*}
  Thus $\mim{\mb{E}} = \shortcentral{\grpu{d}}$ and the first property holds. Now since $\mexp{A}$ and $U^{\otimes n}$ commute for $U \in \grpu{d}$, we get that
  \begin{equation*}
    \mexp{\mexp{A}} = \mexp{\mexp{A}\id{V^{\otimes n}}} = \mexp{A}\mexp{\id{V^{\otimes n}}} = \mexp{A}
  \end{equation*}
  \ie{}, $\mb{E}^2 = \mb{E}$. We first prove the second property in order to see that $\mb{E}$ is indeed an orthogonal projection. This follows immediatly from the fact that the Haar integral commutes with the trace,
  \begin{equation*}
    \begin{split}
      \mtr{\mexp{A}} & = \mtr{\uint{U^{\otimes n}A (U^*)^{\otimes n}}} \\
      & = \uint{\mtr{U^{\otimes n}A(U^*)^{\otimes n}}} \\
      & = \uint{\mtr{A}} = \mtr{A}.
    \end{split}
  \end{equation*}
  Thus we have
  \begin{equation*}
    \begin{split}
      \spl \mexp{A}, \mexp{B} \spr & = \mtr{\mexp{A}^*\uint{U^{\otimes n} B (U^*)^{\otimes n}}} \\
      & = \uint{\mtr{U^{\otimes n}\mexp{A}^* B (U^*)^{\otimes n}}} \\
      & = \uint{\mtr{\mexp{A}^*B}} \\
      & = \mtr{\mexp{A}^*B} = \spl \mexp{A},B\spr.
    \end{split}
  \end{equation*}

  Property 3 is a direct consequence of the fact that all $X,Y \in \shortcentral{\grpu{d}}$ commute with $U^{\otimes n}$ for $U \in \grpu{d}$. 
  
  The fourth property is now a simple calculation
  \begin{equation*}
    \begin{split}
      \mtr{A\mexp{B}} & = \mtr{\mexp{A\mexp{B}}} = \mtr{\mexp{A}\mexp{B}}\\
      &= \mtr{\mexp{\mexp{A}B}} = \mtr{\mexp{A}B} = \mtr{B\mexp{A}}.
    \end{split}
  \end{equation*}

  To prove the fifth property we define for some multiindex $\indi$ the multiset $\set{\indi} = \multiset{i_j}_{j=1}^n$ of all integers appearing as entries of $\indi$. Note that $\indj \neq \sigma(\indi)$ for all $\sigma \in \grps{n}$ if and only if $\set{\indi} \neq \set{\indj}$. Let $i^* \in \set{\indi}\setminus \set{\indj}$,  $\lambda\in\cc$ such that $|\lambda| =1$ and let $t$ be the number of times $i^*$ appears as an entry of $\indi$. The operator
  \begin{equation*}
    T = \lambda \bv{i^*}\dbv{i^*} + \sum_{i\neq i^*}\bv{i}\dbv{i}
  \end{equation*}
  is an element of $\grpu{d}$. By the first property we have
  \begin{equation*}
    \mexp{\bv{\indi}\dbv{\indj}} = \mexp{T^{\otimes n} \bv{\indi}\dbv{\indj}(T^*)^{\otimes n}} = \lambda^t \mexp{\bv{\indi}\dbv{\indj}}.
  \end{equation*}
  For a suitable choice of $\lambda$ this implies that $\mexp{\bv{\indi}\dbv{\indj}} = 0$.
  
  To prove the last property, recall that $S_\rho \in \grpu{d}$. By the first property 
  \begin{equation*}
    \mexp{\bv{\indi}\dbv{\indj}} = \mexp{S_\rho^{\otimes n}\bv{\indi}\dbv{\indj}S_{\rho^{-1}}^{\otimes n}} = \mexp{\bv{\rho[\indi]}\dbv{\rho[\indj]}}.
  \end{equation*}
\end{proof}

\begin{mprop}[Symmetries of the integral]\label{int-sym}
  The integral vanishes unless $\indci = \sigma(\indi)$ and $\indcj = \tau(\indj)$ for some $\sigma, \tau \in \grps{n}$, and has the following symmetries:
  \begin{enumerate}
  \item $\uinti{\indi}{\indj}{\indci}{\indcj} = \uinti{\sigma(\indi)}{\sigma(\indj)}{\tau(\indci)}{\tau(\indcj)}$ ,
  \item $\uinti{\indi}{\indj}{\indci}{\indcj} = \uinti{\pi[\indi]}{\rho[\indj]}{\pi[\indci]}{\rho[\indcj]}$,
  \item $\uinti{\indi}{\indj}{\sigma(\indi)}{\tau(\indj)} = \uinti{\indi}{\indj}{\indi}{\sigma^{-1}(\tau(\indj))}$,
  \item $\uinti{\indi}{\sigma(\indj)}{\indci}{\tau(\indcj)} = \uinti{\sigma^{-1}(\indi)}{\indj}{\tau^{-1}(\indci)}{\indcj}$,
  \item $\uinti{\indi}{\indj}{\indci}{\indcj} = \uinti{\indcj}{\indci}{\indj}{\indi}$,
  \end{enumerate}
  for all $\sigma, \tau \in \grps{n}$ and $\pi , \rho \in \grps{d}$. 
\end{mprop}
\begin{proof}
  By the fifth property of Proposition \ref{prop-e}, we have
  \begin{equation*}
    \uinti{\indi}{\indj}{\indci}{\indcj} = \mtr{\bv{\indci}\dbv{\indi}\mexp{\bv{\indj}\dbv{\indcj}}} = \mtr{\bv{\indj}\dbv{\indcj}\mexp{\bv{\indci}\dbv{\indi}}},
  \end{equation*}
  which is zero unless $\indci = \sigma(\indi)$ and $\indcj = \tau(\indj)$ for some $\sigma, \tau \in \grps{n}$.
  The first symmetry is just a rewording of the fact, that multiplication in $\cc$ is commutative. To prove the second symmetry we apply the fourth and last property of Proposition \ref{prop-e} multiple times,
  \begin{equation*}
    \begin{split}
      \uinti{\indi}{\indj}{\indci}{\indcj} & = \mtr{\bv{\indci}\dbv{\indi}\mexp{\bv{\indj}\dbv{\indcj}}} \\
      & = \mtr{\bv{\indci}\dbv{\indi}\mexp{\bv{\rho[\indj]}\dbv{\rho[\indcj]}}} \\
      & = \mtr{\bv{\rho[\indj]}\dbv{\rho[\indcj]} \mexp{\bv{\indci}\dbv{\indi}}} \\
      & = \mtr{\bv{\rho[\indj]}\dbv{\rho[\indcj]} \mexp{\bv{\pi[\indci]}\dbv{\pi[\indi]}}} \\
      & = \mtr{\bv{\pi[\indci]}\dbv{\pi[\indi]} \mexp{\bv{\rho[\indj]}\dbv{\rho[\indcj]}}} \\
    & = \uinti{\pi[\indi]}{\rho[\indj]}{\pi[\indci]}{\rho[\indcj]}.
    \end{split} 
  \end{equation*}
  The third property is a consequence of the fact, that the representation $\pi$ of $\grps{n}$ commutes with every $U^{\otimes n}$. We get
  \begin{equation*}
    \begin{split}
      \uinti{\indi}{\indj}{\sigma(\indi)}{\tau(\indj)} &= \mtr{\bv{\sigma(\indi)}\dbv{\indi}\mexp{\bv{\indj}\dbv{\tau(\indj)}}} \\
      & = \mtr{\pi(\sigma)\bv\indi \dbv{\indi}\mexp{\bv{\indj}\dbv{\tau(\indj)}}} \\
      & = \mtr{\bv{\indi}\dbv{\indi}\mexp{\bv{\indj}\dbv{\tau(\indj)}\pi(\sigma^{-1})^*}} \\
      & = \mtr{\bv{\indi}\dbv{\indi}\mexp{\bv{\indj}\dbv{\sigma^{-1}(\tau(\indj))}}} \\
      & = \uinti{\indi}{\indj}{\indi}{\sigma^{-1}(\tau(\indj))}.
    \end{split}
  \end{equation*}
  To prove the fourth property note that
  \begin{equation*}
    \begin{split}
      \uinti{\indi}{\sigma(\indj)}{\indci}{\tau(\indcj)} & = \mtr{\bv{\indci}\dbv{\indi}\mexp{\bv{\sigma(\indj)}\dbv{\tau(\indcj)}}} \\
      & = \mtr{\bv{\indci}\dbv{\indi}\mexp{\pi(\sigma)\bv{\indj}\dbv{\indcj} \pi(\tau)^*}} \\
      & = \mtr{\pi(\tau^{-1}) \bv{\indci}\dbv{\indi} \pi(\sigma^{-1})^*\mexp{\bv{\indj}\dbv{\indcj}}} \\
      & = \uinti{\sigma^{-1}(\indi)}{\indj}{\tau^{-1}(\indci)}{\indcj}.
    \end{split}
  \end{equation*}
  The last property is a simple consequence of the fact that the map $U \mapsto U^*$ is itself unitary.
\end{proof}

This allows us to quickly recover some classical results. 

\section{Classical results}
The following section is concerned with the classical theory of the integral \eqref{eq-basic-int}. The results can be found in the book of E. Hewitt and K. A. Ross \cite{hewittross}. 

\begin{mprop}[\cite{hewittross}, p. 116]
  Let $a_j$ and $b_j$ be arbitrary non-negative integers and let $n_j \in \set{1,\dots,d}$ for $j = 1, \dots, d$. Then the integrals
  \begin{enumerate}
  \item $\displaystyle \uint{\prod_{j=1}^d \left(u_{n_jj}\right)^{a_j} \left(\bar{u}_{n_jj}\right)^{b_j}}$,
  \item $\displaystyle \uint{\prod_{j=1}^d \left(u_{jj}\right)^{a_j} \left(\bar{u}_{jj}\right)^{b_j}}$,
  \item $\displaystyle \uint{\prod_{j=1}^d \left(u_{1j}\right)^{a_j} \left(\bar{u}_{1j}\right)^{b_j}}$,
  \end{enumerate}
  vanish unless $a_j = b_j$ for all $j$,
  
\end{mprop}
\begin{proof}
  The first statement follows immediately from the fact that the multiindices
  \begin{equation*}
    \indj = (\underbrace{1,\dots,1}_{a_1}, \dots ,\underbrace{d,\dots,d}_{a_d})
  \end{equation*}
  and
  \begin{equation*}
    \indcj = (\underbrace{1,\dots,1}_{b_1}, \dots ,\underbrace{d,\dots,d}_{b_d})
  \end{equation*}
  have to be permutations of eachother in order for the integral to not vanish. This is only possible if $a_j = b_j$ for every $j$.
  The last two claims are direct consequences of the first. 
\end{proof}

Next we define $\parts{n}{d} = \set{ \mi{a} \in \natnum^{d} \mid \sum_i a_i = n}$, the set of all non-ordered partitions (\textit{compositions}) of $n$ consisting of atmost $d$ parts. 
For $\mi{a} \in \parts{n}{d}$ we define the multinomial coefficient to be
\begin{equation*}
  \binom{n}{\mi{a}} = \frac{n!}{a_1! \cdots a_n!}.
\end{equation*}

For $x_1, \dots, x_n \in \cc$ we write $x^{\mi{a}}$ inplace of $x_1^{a_1} \cdots x_d^{a_d}$. Using this we get that 
\begin{equation*}
  \left(\sum_{i=1}^{d} x_i\right)^n = \sum_{\mi{a} \in \parts{n}{d}}\binom{n}{\mi{a}} x^{\mi{a}}.
\end{equation*}

Furthermore the coefficient of $|x|^{2\mi{a}} = |x_1|^{2a_i} \cdots |x_d|^{2a_d}$ in the sum $\left| \sum_{i=1}^{d}x_i\right|^{2n}$ is given by $\binom{n}{\mi{a}}^2$. To see this note that for any $\mi{a} \in \parts{n}{d}$ the term $|x|^{2\mi{a}}$ occurs exactly once in the double sum
\begin{equation*}
  \begin{split}
    \left|\sum_{i=1}^{d}x_i\right|^{2n}& = \left(\sum_{i=1}^d x_i \right)^n \left(\sum_{i=1}^d \bar{x}_i \right)^n\\
    & = \left(\sum_{\mi{a} \in \parts{n}{d}} \binom{n}{\mi{a}}x^{\mi{a}}\right)^n \left(\sum_{\mi{b} \in \parts{n}{d}} \binom{n}{\mi{b}}\bar{x}^{\mi{b}}\right)^n,
  \end{split} 
\end{equation*}
namely when $\mi{b} = \mi{a}$.

The following two lemmas, although somewhat technical in their nature, are the main ingredients to the computation of the integral $\uinti{\indi}{\indj}{\indci}{\indcj}$ for some special choices of $\indi,\indj,\indci$ and $\indcj$. 

\begin{lemma}[\cite{hewittross}, Lemma 29.7]\label{hr-trig-indep}
  For every $t \in \natnum$ the set of functions
  \begin{equation}\label{trigindp}
    \set{\cos(\theta)^{2t - 2r}\sin(\theta)^{2r} \mid r = 0, \dots, t} 
  \end{equation}
  for $\theta \in [0,\frac{\pi}{2}]$ is linearly independent. 
\end{lemma}

\begin{proof}
  Assume that $\displaystyle \sum_{r = s}^t a_r \cos(\theta)^{2t -2r} \sin(\theta)^{2r} = 0$ for $\theta \in (0,\frac{\pi}{2})$ and that $a_s \neq 0$. Dividing by $\sin(\theta)^{2s}$ yields,
  \begin{equation*}
    \sum_{r = s}^t a_r \cos(\theta)^{2t -2r} \sin(\theta)^{2r-2s} = 0.
  \end{equation*}
  Sending $\theta \rightarrow 0$ implies $a_s = 0$, contradicting our assumption.
\end{proof}

\begin{lemma}[\cite{hewittross}, Lemma 29.8]\label{lemma:hr-tech}
  Let $s$ and $d$ be positive integers, $\mi{a} \in \parts{n}{d}$ and $\mi{n} = \left(n_j\right)_{j=1}^d \in \idns{d}{d}$. If $n_k = n_l$ for some $k,l \in \set{1,\dots,d}$ then we have
  \begin{multline}
    \label{tech-lemma1}
      (a_k + 1) \uint{|u_{n_ll}|^2 \prod_{j=1}^d|u_{n_jj}|^{2a_j}} = \\
      (a_l + 1) \uint{|u_{n_kk}|^2 \prod_{j=1}^d|u_{n_jj}|^{2a_j}}.
    \end{multline}

  \end{lemma}
  \begin{proof}
    We may assume that $k \neq l$, since otherwise the assertion is trivial. Furthermore we may assume without loss of generality that $k=1, l=2$ and $n_l = 1$. To see this, choose permutations $\tau, \sigma \in \grps{d}$ such that $\sigma(n_l) =1$, $\tau(k) = 1$ and $\tau(l) = 2$ and apply the symmetries from Proposition \ref{int-sym}.

    Now, for $\theta$ in $[0,\frac{\pi}{2}]$ define $V_\theta \in \grpu{d}$ by the following rules

    \begin{align*}
      V_\theta e_1 & = \cos(\theta) e_1 + \sin(\theta) e_2, \\
      V_\theta e_2 & = -\sin(\theta)e_1 + \cos(\theta) e_2,\\
      V_\theta e_i & = e_i ~~~\text{for $i \neq 1,2$.}
    \end{align*}
    With this we get
    \begin{multline} \label{hr-tech-lem-eq}
      I := \uint{|u_{11}|^{2a_1+2a_2+2}\prod_{j=3}^d|u_{n_jj}|^{2a_j}} = \\
      \uint{|e_1^*Ue_1|^{2a_1+2a_2+2}\prod_{j=3}^d|e_{n_j}^{*}Ue_j|^{2a_j}} = \\
      \uint{|e_1^*UV_\theta e_1|^{2a_1+2a_2+2}\prod_{j=3}^d|e_{n_j}^{*}UV_\theta e_j|^{2a_j}} = \\
      \uint{|\cos(\theta) u_{11} + \sin(\theta) u_{12}|^{2a_1+2a_2+2}\prod_{j=3}^d|u_{n_jj}|^{2a_j}},
    \end{multline}
    using the convention that empty products are 1.

    Expanding $|\cos(\theta) u_{11} + \sin(\theta) u_{12}|^{2a_1+2a_2+2}$ by the binomial theorem and using Proposition \ref{int-sym} we see that all non-vanashing terms of the resulting sum have integrands of the form (using $t := a_1+ a_2+1$)
    \begin{equation*}
      |\cos(\theta)|^{2r}|\sin(\theta)|^{2t - 2r}\prod_{j=3}^d|u_{n_jj}|^{2a_j},
    \end{equation*}
    for $r = 0, \dots, t$. It follows that the integral in \eqref{hr-tech-lem-eq} is equal to
    \begin{equation}\label{hr-tech-lemma-eq2}
      I = \sum_{r=0}^t \binom{t}{r}^2 \sin(\theta)^{2r}\cos(\theta)^{2t-2r} \uint{|u_{11}|^{2r} |u_{12}|^{2t}\prod_{j=3}^d|u_{n_jj}|^{2a_j}}.
    \end{equation}
    Furthermore we have
    \begin{equation}\label{hr-tech-lemma-eq3}
      I = \left(\sin(\theta)^2 + \cos(\theta)^2\right)^tI = \sum_{r=0}^t \binom{t}{r} \sin(\theta)^{2r}\cos(\theta)^{2t-2r}I.
    \end{equation}
    By Lemma \ref{hr-trig-indep} we can compare the coefficients in the sums \eqref{hr-tech-lemma-eq2} and \eqref{hr-tech-lemma-eq3}.
    Doing this yields
    \begin{equation}
      \label{eq:hr-tech-lemma-coeff1}
      \binom{t}{a_1}^2\uint{|u_{11}|^{2a_1} |u_{12}|^{2a_2}\prod_{j=3}^d|u_{n_jj}|^{2a_j}} = \binom{t}{a_1}I,
    \end{equation}
    for $r=a_1$ and
    \begin{equation}
      \label{eq:hr-tech-lemma-coeff2}
      \binom{t}{a_1+1}^2\uint{|u_{11}|^{2a_1+1} |u_{12}|^{2a_2}\prod_{j=3}^d|u_{n_jj}|^{2a_j}} = \binom{t}{a_1+1}I
    \end{equation}
    for $r = a_1 +1$. Combining \eqref{eq:hr-tech-lemma-coeff1} and \eqref{eq:hr-tech-lemma-coeff2} we get
    \begin{multline*}
      \binom{t}{a_1+1}\uint{|u_{11}|^{2a_1+1} |u_{12}|^{2a_2}\prod_{j=3}^d|u_{n_jj}|^{2a_j}} \\= \binom{t}{a_1} \uint{|u_{11}|^{2a_1} |u_{12}|^{2a_2}\prod_{j=3}^d|u_{n_jj}|^{2a_j}}
    \end{multline*}
    from which the theorem follows. 
  \end{proof}

  \begin{theorem}[\cite{hewittross}, Theorem 29.9]\label{thm:hr-int-form}
    For $\mi{a} \in \parts{n}{d}$ we have
    \begin{equation*}
      \uint{\prod_{j=1}^d|u_{1j}|^{2a_j}} = \frac{(d-1)!}{(n+d-1)!}\prod_{j=1}^da_j!.
    \end{equation*}
  \end{theorem}
  \begin{proof}
    This theorem will be proved inductively. For $n=1$ we already know that $\uint{|u_{1j}|^2} = \frac{1}{d}$. Now let $n \geq 2$ and assume that the theorem holds for all $\mi{a}\in \parts{n-1}{d}$. Using Lemma \ref{lemma:hr-tech} with $l = 1$, $a_l = a_1-1$ and $n_j = 1$ for all $j$ we get
    \begin{equation}
      \label{eq:hr-tech-lemma-app}
      (a_k +1) \uint{\prod_{j=1}^d|u_{1j}|^{2a_j}} = a_1 \uint{|u_{1k}|^{2} |u_{11}|^{2a_1-2} \prod_{j=2}^d|u_{1j}|^{2a_j}}.
    \end{equation}
    Plugging in the induction hypothesis yields
    \begin{equation*}
      \frac{(d-1)!}{(n+d-2)!}(a_1 -1)!\prod_{j=2}^d a_j! = \uint{|u_{11}|^{2a_1 - 2}\prod_{j=2}^d |u_{1j}|^{2a_j}}.
    \end{equation*}
    Since the columns of a unitary matrix have norm $1$, this is equal to
    \begin{multline*}
      \uint{|u_{11}|^{2a_1 - 2}\left(\sum_{k=1}^d |u_{1k}|^2\right)\prod_{j=2}^d |u_{1j}|^{2a_j}} = \\
      \uint{\prod_{j=1}^d |u_{1j}|^{2a_j}} + \sum_{k=2}^d \uint{|u_{1k}|^2 |u_{11}|^{2a_1 - 2}\prod_{j=2}^d |u_{1j}|^{2a_j}}.
    \end{multline*}
    Using \eqref{eq:hr-tech-lemma-app} and rearranging terms we see, that this is the same as
    \begin{equation*}
      \left(1+\sum_{k=2}^d \frac{a_k+1}{a_1}\right)\uint{\prod_{j=1}^d |u_{1j}|^{2a_j}} = \frac{n+d-1}{a_1}\uint{\prod_{j=1}^d |u_{1j}|^{2a_j}}.
    \end{equation*}
    Solving for $\uint{\prod_{j=1}^d |u_{1j}|^{2a_j}}$ finishes the proof. 
  \end{proof}

  \begin{corollary}[\cite{hewittross}, Corollary 29.10]
    For $i,j \in \set{1,\dots,d}$ and $n\in \natnum$ we have
    \begin{equation*}
      \uint{|u_{ij}|^{2n}} = \binom{d+n-1}{d-1}^{-1}.
    \end{equation*}
  \end{corollary}
  \begin{proof}
    This is a direct consequence of Theorem \ref{thm:hr-int-form} using $\mi{a} = (n,0,\dots,0)$ and the symmetries from Proposition \ref{int-sym}.
  \end{proof}

  So far we have only computed special cases by somewhat adhoc methods. Our next goal is to compute $\uinti{\indi}{\indj}{\indci}{\indcj}$ for arbitrary values of $\indi, \indj, \indci$ and $\indcj$.

\section{Weingarten calculus}
In this section we will present the method used by B. Collins and P. Sniady in \cite{collins2006} to derive a formula for the Weingarten function.
Our main tools will be the Schur-Weyl duality (Theorem \ref{thm:swd}), the double centralizer theorem (Theorem \ref{thm:cent}) and the representations $\pi$ and $\rho$ of $\grps{n}$ and $\grpu{d}$ respectively. 

We can use linear extensions of $\pi$ and $\rho$ to map $\gas{n}$ and $\gau{d}$ to subalgebras of $\mend{\vtn}$.
Using the notation from Section \ref{sec:schur-weyl} these subalgebras are denoted by $\mathcal{A} = \pi(\gas{n}) =  \spl \pi(\sigma) \mid \sigma \in \grps{n}\spr$ and $\mathcal{B}= \rho(\gau{d}) = \spl U^{\otimes n} \mid U \in \grpu{d}\spr$. Recall that the group algebra of the symmetric group admits the following decomposition
\begin{equation}
  \label{eq:gau-decomp}
  \begin{split}
    \gas{n} &= \bigoplus_{\lambda \vdash n}\mend{\spv}.
  \end{split}
\end{equation}

There is however no a priori reason why its image $\mathcal{A}$ as a subset of $\mend{\vtn}$ should respect this decomposition. In fact we will see that this does not hold if $d<n$. Not accounting for this was one of the main shortcomings of the initial attempts to compute the integral. 

In order to describe $\mathcal{A}$ we define the following subalgebra of $\gas{n}$
\begin{equation}
  \label{eq:def-cdsn}
  \sncd = \left(\sum_{\fulls} \snmcp \right) \gas{n} = \bigoplus_{\fullss} \mend{\spv}.
\end{equation}
As we will prove shortly, $\pi(\gas{n}) = \mathcal{A}$. The explicit description of $\mathcal{A}$ is what enabled B. Collins and P. Sniady to drop the restriction $d \geq n$ limiting earlier formulas (\eg{} \cite{samuel}). 

Note that by the Schur-Weyl duality the action of $\pi$ can be thought of as a special case of the joint action of $\pi$ and $\rho$, namely $\pi = \pi\times \rho(\cdot,\ids)$. Furthermore by the Schur-Weyl duality we have

\begin{equation*}
  \begin{split}
    \mend{\vtn} & = \mhom{\bigoplus_{\fullss} \spv \otimes \wev}{\bigoplus_{\fullss}\spvm \otimes \wevm} \\
    & \simeq \bigoplus_{\lambda, \mu} \mhom{\spv \otimes \wev}{\spvm \otimes \wevm} \\
    & \simeq \bigoplus_{\lambda, \mu} \left(\spv \otimes \wev\right)^* \otimes \left( \spvm \otimes \wevm \right) \\
    & \simeq \bigoplus_{\lambda, \mu} (\spv)^* \otimes (\wev)^* \otimes \spvm \otimes \wevm \\
    & \simeq \bigoplus_{\lambda, \mu} (\spv)^* \otimes \spvm \otimes (\wev)^* \otimes \wevm \\
    & \simeq \bigoplus_{\lambda, \mu} \mhom{\spv}{\spvm} \otimes \mhom{\wev}{\wevm},
  \end{split}
\end{equation*}

and hence
\begin{equation}\label{eq:pi-equal}
  \begin{split}
    \pi(\gas{n}) &= \pi \times \rho (\gas{n},\ids) = \left(\sum_{\fulls} \spa \otimes \wea(\ids)\right)(\gas{n}) \\
    & = \bigoplus_{\fullss} \spa(\mend{\spv})
  \end{split}
\end{equation}

\begin{lemma}\label{lemma:pi-embeds}
  $\pi$ embeds $\sncd$ into $\mend{\vtn}$, and hence $\pi(\gas{n}) \simeq \sncd $. 
\end{lemma}
\begin{proof}
  Let $x = \sum x_g \delta_g \in \cc[S_n]$. The Schur-Weyl duality implies
  \begin{equation*}
    \pi(x) = \sum_{\fullss}\spa (x)\otimes \id{\wev},
  \end{equation*}
  where $\spa(x) = \sum x_g \spa(g)$.
  This map is essentially the same as the map in Corollary \ref{ga-inj}. The are two minor differences. The range is restricted to components associated to partitions of length at most $d$, and we take tensor products with $\id{\wev}$, but neither of those operations impairs injectivity. 

  Since $\pi$ is injective, each $\spa: \mend{\spv} \fa \gas{n}$ has to be injective. Hence equation (\ref{eq:pi-equal}) implies  that $\pi(\gas{n}) = \sncd$.
\end{proof}

This directly leads to our next proposition.
Recall that for an inclusion of algebras $\mathcal{M} \subseteq \mathcal{N}$ a \textit{conditional expectation} is a $\mathcal{M}$-bimodule map $\mb{E}: \mathcal{N} \fa \mathcal{M}$ such that $\mexp{\id{\mathcal{N}}} = \id{\mathcal{M}}$.
\begin{mprop}
  $\mb{E}$ is a conditional expectation of $\mend{\vtn}$ onto $\sncd$.
\end{mprop}
\begin{proof}
  By the preceeding discussion we get that $\mathcal{A} = \sncd$. By theorem \ref{thm:cent} and Proposition \ref{prop-e} we get $\sncd = \mathcal{A} = \shortcentral{\grpu{d}} = \mim{\mb{E}}$. Since $\mexp{\id{\vtn}} = \id{\vtn}$ property three of Proposition \ref{prop-e} implies that $\mb{E}$ is a conditional expectation of $\mend{\vtn}$ onto $\sncd$.
  
\end{proof}

Next we can define a special element of $\gas{n}$ which will allow us to derive an explicit formula for the Weingarten function. For $A \in \mend{\vtn}$ we define
\begin{equation}
  \label{eq:def-phi}
  \Phi(A) = \sum_{\tau \in\grps{n}} \mtr{A\pi(\tau^{-1})} \delta_\tau.
\end{equation}
\begin{lemma}
  $\gagt{A}$ is a $\gas{n}-\gas{n}$ bimodule homomorphism, in the sense that
  \begin{equation*}
    \begin{split}
      & \gagt{A\pi(\sigma)} = \gagt{A}\cdot \sigma, \\
      & \gagt{\pi(\sigma)A} = \sigma \cdot \gagt{A}.
    \end{split}
  \end{equation*}
\end{lemma}
\begin{proof}
  Note that
  \begin{equation*}
    \begin{split}
      \gagt{A\pi(\sigma)} &= \sum_{g \in S_n}\mtr{A\pi(g^{-1}\sigma)}\delta_g\\
      & = \sum_{g \in S_n}\mtr{A\pi(g^{-1})}\delta_{g\sigma}\\
      & = \gagt{A}\cdot \sigma
    \end{split}
  \end{equation*}
  In the exact same way we can show that
  \begin{equation*}
    \gagt{\pi(\sigma)A} = \sigma \cdot \gagt{A}.
  \end{equation*}
\end{proof}

\begin{mprop}[Weingarten function]
  We have $\gagt{\ids} = \charpi$. Furthermore $\charpi$ is an invertible element of $\gas{n}$ and its inverse is given by
  \begin{equation}
    \label{eq:def-wein}
    \wein = \frac{1}{(n!)^2}\sum_{\fullss}\frac{\mdim{\spv}^2}{\mdim{\wev}} \charspecht,
  \end{equation}
  where $\charspecht$ denotes the character of the irreducible representation $\spv$. 
\end{mprop}
\begin{remark}
  As we will see shortly $\wein$ is the Weingarten function and will allow us to compute the integral $\uinti{\indi}{\indj}{\indci}{\indcj}$. Note that the Weingarten function is constant on conjugacy classes, as a sum of class functions. This means that the Weingarten function essentially is a function of partitions of $n$.  
\end{remark}

\begin{proof}
  Note that
  
  \begin{equation*}
    \gagt{\ids} = \sum_{\tau\in\grps{n}}\mtr{\pi(\tau^{-1})}\delta_\tau = \sum_{\tau \in \grps{n}} \charpiarg{\tau^{-1}}\delta_\tau = \charpi,
  \end{equation*}
  where the last equality is due to the fact that $\tau$ and $\tau^{-1}$ are conjugate.

  Since $\charpi(\tau) = \charjointarg{\tau}{\ids}$ the Schur-Weyl duality implies that
  \begin{equation*}
    \begin{split}
      \charpi(\tau) & = \sum_{\fullss} \charjointdecomparg{\tau}{\ids} \\
      & = \sum_{\fullss} \charspechtarg{\tau} \charweylarg{\ids} \\
      & = \sum_{\fullss} \charspechtarg{\tau} \mdim{\wev} \\
      & = n!\sum_{\fullss}\frac{\mdim{\wev}}{\mdim{\spv}} \snmcp(\tau)
    \end{split}
  \end{equation*}
  Since this sum is direct and $\snmcp$ acts as the identity on the $\lambda$-isotypic component of $\sncd$, all we have to do in order to invert $\charpi$ is taking the reciprocal of the coefficient of $\snmcp$. This yields
  \begin{equation*}
    \text{Wg}: = \charpi^{-1} = \frac{1}{(n!)^2}\sum_{\fullss}\frac{\mdim{\spv}^2}{\mdim{\wev}}\charspecht.
  \end{equation*}
\end{proof}

The following identities will allow us to derive an explicit formula for the integral $\uinti{\indi}{\indj}{\indci}{\indcj}$.
From here on we identify elements of $\sncd$ with elements of $\mathcal{A} = \pi(\sncd) = \mim{\mb{E}}$ and vice versa.  
Since $\pi: \sncd \fa \mend{\vtn}$ is injective this is justified.

\begin{mprop}\label{prop:phi-prop}
  We have
  \begin{enumerate}
  \item $\gagt{A} = \mb{E}(A)\gagt{id{}}$.
  \item $\mim{\Phi} = \sncd$.
  \item $\Phi(A\mb{E}(B)) = \Phi(A) \Phi(B) \gagt{\id{}}^{-1}$.
  \end{enumerate}
\end{mprop}

\begin{proof}
  Since $\mtr{\mb{E}(A)} = \mtr{A}$ and $\mb{E}(A) \in \sncd$, we have
  \begin{equation*}
    \begin{split}
      \gagt{A} & = \sum_{\tau \in S_n}\mtr{A\pi(\tau^{-1})}\delta_\tau = \sum_{\tau \in S_n}\mtr{\mb{E}(A)\pi(\tau^{-1})}\delta_\tau \\
      & = \gagt{\mb{E}(A)} = \gagt{\mb{E}(A)\ids} = \mb{E}(A)\gagt{\ids}.
    \end{split}
  \end{equation*}
  This directly implies that $\mim{\Phi} = \sncd$, because \gagt{\id{}} is invertible and $\mim{\mb{E}} = \sncd$.

  The third equality is a direct consequence of the first one.
  \begin{equation*}
    \begin{split}
      \gagt{A\mb{E}(B)} &= \mb{E}(A\mb{E}(B))\gagt{\ids} = \mb{E}(A)\mb{E}(B) \gagt{\ids} \\
      & = \gagt{A}\gagt{B}\gagt{\ids}^{-2}\gagt{\ids} \\
      & = \gagt{A}\gagt{B}\gagt{\ids}^{-1}. 
    \end{split}
  \end{equation*}

\end{proof}

Now we are finally able to compute the integral.

\begin{theorem}[Collins \& Sniady, \cite{collins2006}]\label{thm:uint-colsn}
  For $\indi, \indj, \indci, \indcj \in \idn$ we have
  
  \begin{equation}
    \label{eq:int-formula}
    \uinti{\indi}{\indj}{\indci}{\indcj} = \sum_{\sigma,\tau \in \grps{n}} \delta_{\sigma(\indi),\indci} \delta_{\tau(\indj),\indcj}\weinarg{\tau\sigma^{-1}}.
  \end{equation}
\end{theorem}
\begin{proof}
  Recall that
  \begin{equation*}
    \uinti{\indi}{\indj}{\indci}{\indcj}= \dbv{\indi}\mexp{\bv{\indj}\dbv{\indcj}}\bv{\indci} = \mtr{\bv{\indci}\dbv{\indi}\mexp{\bv{\indj}\dbv{\indcj}}}. 
  \end{equation*}
  If we define $A = \bv{\indci}\dbv{\indi}$ and $B = \bv{\indj}\dbv{\indcj}$ we get that the integral is equal to $\gagt{A\mexp{B}}(e)$ \ie{}, to the coefficient of $\delta_e$ in $\gagt{A\mexp{B}} = \gagt{A}\gagt{B}\gagt{\ids}^{-1}$.

  Note that
  \begin{equation*}
    \begin{split}
      \gagt{A} & = \sum_{\sigma \in \grps{n}} \mtr{\bv{\indci}\dbv{\indi}\pi(\sigma^{-1})} \delta_\sigma\\
      & = \sum_{\sigma\in\grps{n}} \dbv{\indi}\pi(\sigma^{-1})\bv{\indci} ~ \delta_\sigma \\
      & = \sum_{\sigma \in \grps{n}} \delta_{\indi,\sigma^{-1}(\indci)} ~ \delta_\sigma \\
      & = \sum_{\sigma \in \grps{n}} \delta_{\sigma(\indi),\indci} ~ \delta_\sigma
    \end{split}
  \end{equation*}
  Similarly
  \begin{equation*}
    \gagt{B} = \sum_{\tau \in \grps{n}} \delta_{\indcj,\tau^{-1}(\indj)}~ \delta_\tau= \sum_{\tau \in \grps{n}} \delta_{\indcj,\tau(\indj)}~ \delta_{\tau^{-1}}.
  \end{equation*}
  Putting everything together we get
  \begin{equation*}
    \begin{split}
      \gagt{A}\gagt{B}\wein &= \sum_{\sigma \in \grps{n}} \delta_{\sigma(\indi),\indci} ~ \delta_\sigma \sum_{\tau \in \grps{n}} \delta_{\tau(\indj),\indcj}~ \delta_{\tau^{-1}} \sum_{x\in \grps{n}} \weinarg{x} \delta_x \\
      &= \sum_{\sigma,\tau,x \in \grps{n}} \delta_{\sigma(\indi),\indci}\delta_{\tau(\indj),\indcj}\weinarg{x} ~ \delta_{\sigma \tau^{-1} x}. 
    \end{split}
  \end{equation*}
  Since $\sigma \tau^{-1} x = e$ if and only if $x = \tau \sigma^{-1}$ we get that the coefficient of $\delta_e$ is equal to
  \begin{equation*}
    \sum_{\sigma,\tau \in \grps{n}} \delta_{\sigma(\indi),\indci} \delta_{\tau(\indj),\indcj}\weinarg{\tau\sigma^{-1}}.
  \end{equation*}
  
\end{proof}

\section{Jucys-Murphy elements and the Moore-Penrose inverse}

Next we discuss the connection between the Weingarten function and the Jucys-Murphy elements. This result is due to J. I. Novak \cite{novak2009}. Recall that the Weingarten function $\wein$ is the inverse of the character $\charpi$. Since characters are constant on conjugacy classes, Lemma \ref{lemma:charpi} implies that $\charpi$ can be written as
\begin{equation}
  \label{eq:charpi}
  \charpi = \wein^{-1}  = \sum_{\mu \in \youngpos{n}} d^{l(\mu)} C_\mu.
\end{equation}

\begin{theorem}[Novak, \cite{novak2009}]\label{thm:jucys-novak}
  For $d \geq n$ the Weingarten function equals
  \begin{equation}
    \label{eq:novak-wein}
    \wein = (d+J_1)^{-1}\cdots (d+J_n)^{-1},
  \end{equation}
  where $J_k$ denotes the $k$-th Jucys-Murphy element. 
\end{theorem}
\begin{proof}
  From the discussion preceding Lemma \ref{lemma:jucys-ev} we get that the eigenvalues of $J_k$ (and hence of $-J_k$) all lie in the inverval $[-n+1,n-1]$.
  Hence for $d \geq n$ the element $d+J_k$ is invertible by a Neumann series.
  
  In view of equation \eqref{eq:charpi} and Jucys' theorem (Theorem \ref{thm:jucys}) we write
  \begin{equation*}
    \begin{split}
      (d+J_1)\cdots (d+J_n) &=  \sum_{k=0}^n d^{n-k} e_k(J_1,\dots,J_n) \\
      &= \sum_{k=0}^n \sum_{\substack{\mu \in \youngpos{n} \\l(\mu) = n-k}} d^{n-k}C_\mu \\
      & = \sum_{\mu \in \youngpos{n}}d^{l(\mu)}C_\mu \\
      &= \charpi = \wein^{-1}.
    \end{split}
  \end{equation*}
\end{proof}

The following self-contained, inductive proof of the equality $\charpi = (d+J_1) \cdots (d+J_n)$ is due to P. Zinn-Justin \cite{zinnjustin} and completely avoids the usage of Theorem \ref{thm:jucys}.

\begin{proof}[Alternative proof of theorem \ref{thm:jucys-novak}]
  Note that by expanding the product we get a sum with $n!$ terms. We will inductively establish a one-to-one correspondense between the terms of both sides of the equation $\charpi = (d+J_1) \cdots (d+J_n)$.

  In the base case $n =1$ the equation reduces to $d = d$.
  For $n\geq 2$ and  $\sigma \in \grps{n}$ we distinguish two cases:
  \begin{enumerate}
  \item $\sigma(n) = n$.
  \item $\sigma(n) \neq n$. 
  \end{enumerate}

  In the first case we use the induction hypothesis to see that the term corresponding to $\restr{\sigma}{\set{1,\dots,n-1}}$ is among the terms of the expansion of $(d+J_1) \cdots (d+J_{n-1})$. Since $\sigma$ has one more cycle than $\restr{\sigma}{\set{1,\dots,n-1}}$ we multply it with the $d$ in the remaining term $d + J_n$.
  
  In the second case the decomposition of $\sigma$ into disjoint cycles is of the following form
  \begin{equation*}
    \sigma = \sigma_1 \cdots\sigma_k \cdot(\sigma(n), n , \sigma^{-1}(n), \dots).
  \end{equation*}
  Now we define the permutation $\tau$ to be
  \begin{equation*}
    \tau := \sigma \cdot (n, \sigma^{-1}(n)) = \sigma_1 \cdots \sigma_k \cdot(\sigma(n),\sigma^{-1}(n), \dots) (n),
  \end{equation*}
and apply the induction hypothesis to $\restr{\tau}{\set{1,\dots,n-1}}$. Since $\restr{\tau}{\set{1,\dots,n-1}}$ and $\sigma$ have the same number of cycles  we get $d^{\#(\sigma)} \sigma$ by multiplying the term corresponding to $\restr{\tau}{\set{1,\dots,n-1}}$ with the transposition $(\sigma^{-1}(n),n)$ in $J_n$.  
\end{proof}
Yet another approach to the computation of the Weingarten function is due to P. Zinn-Justin \cite{zinnjustin}. Since we already know that the Weingarten function is the inverse of the character $\charpi$, we can try the following \textit{ansatz}. We define the Gram matrix $G$ by
\begin{equation}
  \label{eq:def-gram}
  G_{\sigma \tau} =  \spl  \pi(\sigma), \pi(\tau) \spr  = \mtr{\pi(\sigma)^*\pi(\tau)} = d^{\#(\sigma^{-1} \tau)}.  
\end{equation}
If $G$ is invertible,  the first row (and hence also the first column) of its inverse $W = G^{-1}$ would contain all the values of the Weingarten function. If $G$ is not invertible we can still define $W$ to be the Moore-Penrose inverse of $G$ \ie{}, the unique matrix satisfying the following properties
\begin{equation}
  \label{eq:weinmatrix-mpi}
  \begin{split}
    GWG &= G, \\
    WGW &= W, \\
    (WG)^* &= WG, \\
    (GW)^* & = GW.
  \end{split}
\end{equation}
The existence of the Moore-Penrose inverse can be directly verified by considering either the singular value decomposition, or a rank factorization of $G$. 
Note that if $G$ is invertible, the inverse and the Moore-Penrose inverse agree.
The matrix $W$ has been called Weingarten matrix by P. Zinn-Justin.
As we will see the name is justified and the definition of $W$ is perfectly consistent with Theorem \ref{thm:uint-colsn} and the theory outlined so far.

Note that $\charpi$ can be written as
\begin{equation*}
  \charpi = \sum_{\sigma \in \grps{n}} d^{\#\sigma} \delta_\sigma.
\end{equation*}
This looks strikingly similar to the Gram matrix $G$ just defined. In fact, if we let $\charpi$ act on $\gas{n}$ by left (or right multiplication), we get $L(\charpi) = G$, with respect to standard basis $\grps{n}$ of $\gas{n}$.

We can now use the properties of Young's orthogonal idempotents $e_T$ outlined in Section \ref{section:reptheo} to simplify this. 
Note that the equation
\begin{equation*}
  \charpi = \prod_{k=1}^n (d+J_k)
\end{equation*}
remains true even if $d<n$. Multiplying by $1 = \sum_{T:|T| = n} e_T$ we get
\begin{equation*}
  \begin{split}
    \prod_{k=1}^n 1(d+J_k) & = \prod_{k=1}^n \sum_{T:|T| = n} e_T(d+J_k) \\
    & = \prod_{k=1}^n \sum_{T:|T| = n } e_T(d+ \cont{T_k}).
  \end{split}
\end{equation*}
Ordering the tableaux by the shape of their underlying Young diagram and plugging in the definition of $\cont{T_k}$ we get
\begin{equation}
  \label{eq:zj-formula}
  \charpi  =  \prod_{k=1}^n \sum_{T:|T| = n } e_T(d+ \cont{T_k})  = \sum_{\lambda \vdash n} c_\lambda \snmcp, 
\end{equation}
with $c_\lambda = \prod_{(i,j) \in \lambda } (d+j-i)$.
If $G$ is invertible we immediately see that its inverse equals
\begin{equation*}
  G^{-1} = \sum_{\lambda \vdash n} c_\lambda^{-1} L(\snmcp). 
\end{equation*}
If $G$ is not invertible we can still define $w$ 
\begin{equation}\label{eq:weinmatrix-mp}
  w = \sum_{\substack{\lambda \vdash n \\ c_\lambda \neq 0}} c_\lambda^{-1} \snmcp,
\end{equation}
and consider $W = L(w)$. 
Using
\begin{equation*}
  w\charpi = \sum_{\substack{\lambda \vdash n \\ c_\lambda \neq 0}}  \snmcp
\end{equation*}
we get $WGW = W$ which in turn implies that $W$ satisfies the conditions of the Moore-Penrose inverse. Note that for $d \geq n$ all the $c_\lambda$ are non-zero since the difference $j-i$ can never be smaller than $-n+1$. On the other hand if $d < n $ one can certainly find a partition $\lambda \vdash n$ such that the Young diagram of $\lambda$ contains a box $(i,j)$ where $d + j -i = 0$. Hence $d \geq n$ is a necessary and sufficient condition for $G$ to be invertible. Furthermore the condition $l(\lambda) \leq d$ in the definition of the Weingarten function \eqref{eq:def-wein} is equivalent to $c_\lambda \neq 0$.

Using Theorem \ref{thm:hook-length} we write
\begin{equation}
  \label{eq:cl-equiv}
  \begin{split}
    \mdim{\wev} &= \prod_{(i,j) \in \lambda} \frac{d+j-i}{h_\lambda(i,j)} \\
    & = \prod_{(i,j) \in \lambda} \frac{1}{h_\lambda(i,j)} c_\lambda = \frac{\mdim{\spv}}{n!} c_\lambda.
  \end{split}
\end{equation}
Since the minimal central character $\snmcp$ can be written as $\snmcp = \frac{\mdim{\spv}}{n!}\charspecht$, this yields
\begin{equation*}
  \begin{split}
    w &= \sum_{\substack{\lambda \vdash n \\ c_\lambda \neq 0}} \frac{1}{c_\lambda} \snmcp = \frac{1}{(n!)^2}\sum_{\fullss} \frac{\mdim{\spv}^2}{\mdim{\wev}}\charspecht = \wein.
  \end{split}
\end{equation*}

Since the Weingarten function is constant on conjugacy classes, the matrix $W$ is symmetric.  Hence it is also the matrix of the right-regular action of $\wein$. 
Summing up, we get the following theorem.

\begin{theorem}[Zinn-Justin, \cite{zinnjustin}]
  The Moore-Penrose inverse $W$ of the Gram matrix $G$ corresponds to the matrix of the regular action of the Weingarten function on $\gas{n}$ \ie{}, $W_{\tau\sigma} = \weinarg{\tau^{-1}\sigma}$. 
\end{theorem}

Given the form of $W$ in \eqref{eq:weinmatrix-mp} and the fact that $\snmcp$ acts as the identity on the $\lambda$-isotypical components of $\gas{n}$, we can immediately describe the eigenstructure of $W$.

\begin{corollary}[Eigenstructure of $W$]
  $W$ has $\mdim{\spv}^2$ eigenvectors with eigenvalue $c_\lambda^{-1}$ for $l(\lambda) \leq d$, and the nullity of $W$ is given by
  \begin{equation*}
    \nullity{W} = \sum_{\substack{\lambda \vdash n \\ l(\lambda) > d}} \mdim{\spv}^2.
  \end{equation*}
\end{corollary}

\begin{mprop}[B. Collins, S. Matsumoto, N. Saad, \cite{colmatsaad}] 
  \begin{equation}
    \label{eq:sum-wein}
    \sum_{\sigma \in \grps{n}} \weinarg{\sigma} = \frac{1}{d(d+1) \cdots (d+n-1)}.
  \end{equation}
\end{mprop}
\begin{proof}
  On the one hand this is a direct consequence of Theorem \ref{thm:hr-int-form} and Theorem \ref{thm:uint-colsn} (setting $a_j = 1$ for all $j$).
  On the other hand we can apply the trivial representation $\trivnoarg: \gas{n} \fa \cc$, defined by $\sigma \mapsto 1$, to the identity
  \begin{equation*}
    \charpi \wein \charpi = \charpi.
  \end{equation*}
  Since
  \begin{equation*}
   \charpi = \prod_{k=1}^n(d + J_k) ~~ \text{and} ~~~  \triv{\prod_{k=1}^n(d + J_k)} = \prod_{k=1}^n(d + k-1)
  \end{equation*}
  the theorem follows.  
\end{proof}

The Weingarten formula can also be used to compute the joint expectation of the entries of a complex hermitian random matrix invariant under unitary conjugation.
To this end we introduce a new notation. Let $W$ be a $d\times d$ complex hermitian random variable such that $UWU^*$ has the same distribution as $W$ for all $U \in \grpu{d}$. For $\tau \in \grps{n}$ with cycle type $\lambda(\tau) = (\tau_1, \dots, \tau_{l})$ we define
\begin{equation}
  \label{eq:def-tr-sn}
  \mtrsn{\tau}{W} = \mtr{W^{\tau_1}} \cdots \mtr{W^{\tau_l}}.
\end{equation}

\begin{theorem}[Collins, Matsumoto, Saad, \cite{colmatsaad}]\label{thm:inv-her-rm}
  Let $W$ be as above. For $\indi, \indj \in \idn$ we have
  \begin{equation}
    \label{eq:inv-her-rm}
    \begin{split}
      \mexp{W_{\indi\indj}^{\otimes n}} &= \mexp{W_{i_1j_1} \cdots W_{i_nj_n}} \\
      & = \sum_{\sigma, \tau \in \grps{n}}  \delta_{\indi, \sigma(\indj)} \weinarg{\sigma^{-1}\tau} \mexp{\mtrsn{\tau}{W}}.
    \end{split}
  \end{equation}
\end{theorem}

The proof hinges on the so called \textit{spectral theorem for unitary invariant random matrices}.
For a proof  of this theorem we refer to B. Collins and C. Male \cite{collinsmale}.

\begin{theorem}\label{thm:spec-u-inv}
  Let $W$ be a $d\times d$ complex hermitian or unitary matrix whose distribution is invariant under unitary conjugation.
  Then $W = UDU^*$,  where
  \begin{enumerate}
  \item $U$ is a unitary Haar random Matrix,
  \item $D$ is a diagonal matrix containing the eigenvalues of $W$ in increasing order if $W$ is hermitian, or in increasing order of the argument $\theta \in [-\pi,\pi)$ if $W$ is unitary,
  \item $U$ and $D$ are independent. 
  \end{enumerate}
\end{theorem}

\begin{proof}[Proof (Theorem \ref{thm:inv-her-rm})]
  Let $W = UDU^*$ as in Theorem \ref{thm:spec-u-inv} and denote the eigenvalues of $W$ with $c_1, \dots, c_d$ \ie, $D = \mdiag{c_1,\dots,c_d}$.
  Using the independence of $U$ and $D$, and that $D^{\otimes n}$ can be written as
  \begin{equation*}
    D^{\otimes n} = \sum_{\mi{r} \in \idn} c_{r_1} \cdots c_{r_n} \bv{\mi{r}}\dbv{\mi{r}}
  \end{equation*}
  we get
  \begin{equation*}
    \begin{split}
      \mexp{W_{i_1j_i}\cdots W_{i_nj_n}} &= \sum_{\mi{r} \in \idn} \mexp{c_{r_1} \cdots c_{r_n}} \uinti{\indi}{\mi{r}}{\indj}{\mi{r}} \\
      & = \sum_{\sigma, \tau \in \grps{n}} \delta_{\indi, \sigma(\indj)} \weinarg{\sigma^{-1}\tau} \sum_{\mi{r}\in \idn} \delta_{\mi{r},\tau(\mi{r})} \mexp{c_{r_1} \cdots c_{r_n}} \\
      & = \sum_{\sigma, \tau \in \grps{n}} \delta_{\indi, \sigma(\indj)} \weinarg{\sigma^{-1}\tau} \mexp{\mtrsn{\tau}{W}}.
    \end{split}
  \end{equation*}
\end{proof}

\section{A new result}

In \cite{diaconis} P. Diaconis and S. N. Evans computed among other things
\begin{equation*}
  \uint{|\mtr{U^k}|^{2n}} = k^n n!
\end{equation*}
for $1 \leq kn \leq d$. In his 2014 review article \cite{linzhang} L. Zhang mentions, that the case $kn > d$ is still open. At least for $k=1$ we are able to generalize this. 

\graffito{The techniques used in \cite{diaconis} can also be used to prove Proposition \ref{prop:tr}, however the proof outlined here is consistent with the theory developed by B. Collins and P. Sniady.}
\begin{mprop}\label{prop:tr}
  \begin{equation*}
    \uint{|\mtr{U}|^{2n}}  = \sum_{\fullss}\mdim{\spv}^2.
  \end{equation*}
\end{mprop}
\begin{proof}
  We start by noting that
  \begin{equation*}
    \begin{split}
      \uint{|\mtr{U}|^{2n}} &= \uint{\left(\sum_{k,l = 1}^d u_{kk}\bar{u}_{ll}\right)^n} \\
      & = \sum_{\indi,\indj} \dbv{\indi}\mexp{\bv{\indi}\dbv{\indj}}\bv{\indj}.
    \end{split}
  \end{equation*}
  
  Since $\mathbb{E}$ is the orthogonal projection onto $\shortcentral{\grpu{d}}$, the dimension of $\shortcentral{\grpu{d}}$ is equal to the trace of $\mathbb{E}$,
  \begin{equation*}
    \begin{split}
      \mdim{\shortcentral{\grpu{d}}} & = \sum_{\indi,\indj} \spl \bv{\indi}\dbv{\indj}, \mexp{\bv{\indi}\dbv{\indj}} \spr \\
      & = \sum_{\indi, \indj} \mtr{(\bv{\indi}\dbv{\indj})^*\mexp{\bv{\indi}\dbv{\indj}}} \\
      & = \sum_{\indi,\indj} \dbv{\indi}\mexp{\bv{\indi}\dbv{\indj}}\bv{\indj}.
    \end{split}
  \end{equation*}

  Using Lemma \ref{lemma:pi-embeds} and Theorem \ref{thm:cent} we get that
  \begin{equation*}
    \mdim{\shortcentral{\grpu{d}}} = \mdim{\sncd} = \sum_{\fullss}\mdim{\spv}^2.
  \end{equation*}
\end{proof}


The integral $\uinti{\indi}{\indi}{\indj}{\indj}$ vanishes unless $\indj = \sigma(\indi)$ for some $\sigma \in \grps{n}$, in which case the integral is equal to $\uinti{\indi}{\indi}{\indi}{\indi}$ by Proposition \ref{int-sym}. Thus for every $\indi$, $\indj$ only takes values in $\orb{\indi}$. Since $\stab{\indi}$ is a subgroup of $\grps{n}$ we have
\begin{equation*}
  \begin{split}
    \uinti{\indi}{\indi}{\indi}{\indi} &= \sum_{\sigma, \tau \in \grps{n}} \delta_{\indi,\sigma(\indi)}  \delta_{\indi,\tau(\indi)} \weinarg{\sigma\tau^{-1}}  \\
    & = |\stab{\indi}| \sum_{\sigma \in \stab{\indi}} \weinarg{\sigma}.
  \end{split}
\end{equation*}

The Orbit-Stabilizer theorem yields
\begin{equation}\label{eq:trace-stab-sum}
  \begin{split}
    \uint{|\mtr{U}|^{2n}} &= \sum_{\indi,\indj} \uinti{\indi}{\indi}{\indj}{\indj}\\
    &=  \sum_{\indi} |\orb{\indi}|\, \uinti{\indi}{\indi}{\indi}{\indi} \\
    &= \sum_{\indi} |\orb{\indi}|\,|\stab{\indi}| \sum_{\sigma \in \stab{\indi}} \weinarg{\sigma}\\
    &= n! \sum_{\indi} \sum_{\sigma \in \stab{\indi}} \weinarg{\sigma}.
  \end{split}
\end{equation}

Hence $\sum_{\indi} \sum_{\sigma \in \stab{\indi}} \weinarg{\sigma}$ is a measure of how many partitions we lose if we restrict the length of the partitions to be at most $d$. 





\begin{remark}
  Proposition \ref{prop:tr} does not readily generalize to the case $k > 1$. If we denote with $u^k_{ij}$ the $(i,j)$-th entry of $U^k$ we get
  \begin{equation}
    \label{eq:k>1}
    \begin{split}
      \uint{|\mtr{U^k}|^{2n}} &=  \uint{\left(\sum_{r,s = 1}^d u_{rr}^k\bar{u}_{ss}^k\right)^n} \\
      & = \sum_{\indi,\indj} \uint{u_{i_1i_1}^k \cdots u_{i_ni_n}^k\bar{u}_{j_1j_1}^k \cdots \bar{u}_{j_nj_n}^k}
    \end{split}
  \end{equation}
  Plugging in
  \begin{equation*}
    u^k_{ij} = \sum_{p_1, \dots, p_{k-1} = 1}^d u_{i_1p_1} \cdots u_{p_{k-1}i_1}
  \end{equation*}
  we get
  \begin{equation*}
    \begin{split}
      \sum_{\indi,\indj} \uint{&\left(\sum_{p_1, \dots, p_{k-1}=1}^d u_{i_1p_1} \cdots u_{p_{k-1}i_1}  \right) \cdots \left(\sum_{s_1, \dots, s_{k-1}=1}^d u_{i_ns_1} \cdots u_{s_{k-1}i_n}  \right)\\
      &\left(\sum_{q_1, \dots, q_{k-1}=1}^d \bar{u}_{j_1q_1} \cdots \bar{u}_{q_{k-1}j_1}  \right) \cdots \left(\sum_{r_1, \dots, r_{k-1}=1}^d \bar{u}_{j_nr_1} \cdots \bar{u}_{r_{k-1}j_n}  \right)}.
    \end{split}
  \end{equation*}
  Performing the multiplications in the integral we see that
  \begin{equation*}
    \uint{|\mtr{U^k}|^{2n}} = \sum \uinti{\indi}{\indj}{\indci}{\indcj}
  \end{equation*}
  where $\indi, \indj, \indci$ and $\indcj$ are multiindices in $\idns{d}{nk}$ of the form
  \begin{equation*}
    \begin{split}
      \indi  &= (i_1, p_1,p_2, \dots,p_{k-1}, \dots , i_n, s_1 ,s_2, \dots, s_{k-1}),\\
      \indj  &= (p_1,p_2, \dots, p_{k-1},i_1, \dots , s_1, s_2,  \dots, s_{k-1}, i_n), \\
      \indci &= (q_1, q_2, \dots, q_{k-1}, j_1, \dots, r_1, r_2, \dots, r_{k-1}, j_n), \\
      \indcj &= (j_1, q_1, q_2, \dots, q_{k-1}, \dots, j_n, r_1, r_2, \dots, r_{k-1}). 
    \end{split}
  \end{equation*}

  Note that $\indj = \gamma(\indi)$ and $\indcj  = \gamma^{-1}(\indci)$, where $\gamma = \gamma_1\dots \gamma_k \in \grps{nk}$ and $\gamma_l$ is the cyclic permutation given by
  \begin{equation}
    \label{eq:def-beta}
    \gamma_l = (ln, ln -1, \dots , ln - k+1).
  \end{equation}
  Intuitively speaking, we split the numbers $1, \dots, kn$ into $n$ blocks of length $k$ and permute the $l$-th block cyclically according to $\gamma_l$.  
  This yields
  \begin{equation*}
     \uint{|\mtr{U^k}|^{2n}} = \sum_{\indi, \indj \in \idns{d}{kn}} \uinti{\indi}{\gamma(\indi)}{\gamma(\indj)}{\indj}.
   \end{equation*}
   Hence the the simplifications performed in the proof of Proposition \ref{prop:tr} do not carry over.
\end{remark}